\documentclass[12pt]{article}
\usepackage{amsfonts,amssymb,amsmath}
\usepackage{epsfig}
\usepackage{color}
\newtheorem{theorem}{Theorem}[section]
\newtheorem{lemma}[theorem]{Lemma}
\newtheorem{definition}{Definition}[section]
\newcommand{\R}{{\mathbb R}}
\newcommand{\C}{{\mathbb C}}
\newcommand{\Z}{{\mathbb Z}}
\newcommand{\N}{{\mathbb N}}
\newcommand{\be}{\begin{equation}}
\newcommand{\ee}{\end{equation}}
\def\weakly{\rightharpoonup}
\def\weakstar{\stackrel{\ast}{\rightharpoonup}}
\def\wt{\widetilde}
\def\eps{\epsilon}
\def\Psieps{\Psi_{\eps}}
\def\Heps{H_\eps}
\def\Weps{W_{\eps}}
\def\Feps{F_{\eps}}
\def\calF{\mathcal{F}}
\def\calD{\mathcal{D}}
\def\calS{\mathcal{S}}
\def\calA{\mathcal{A}}
\def\calM{\mathcal{M}}
\def\L{\mathrm{L}}
\def\H{\mathrm{H}}
\def\W{\mathrm{W}}
\def\Ce{\mathrm{C}}
\def\d{\,\mathrm{d}}
\def\dd{\mathrm{d}}
\def\dist{\mathrm{dist}}
\def\supp{\mathrm{supp}\,}
\renewcommand{\div}{\mathrm{div}}
\renewcommand{\Re}{\mathrm{Re}}
\newcommand{\ds}{\displaystyle}
\begin{document}
\begin{center}
{\LARGE \bf Passage from} \\ 
{\LARGE \bf quantum to classical molecular dynamics} \\
{\LARGE \bf in the presence of Coulomb interactions} 
\\[5mm]
{\sc 
Luigi Ambrosio\footnote{Scuola Normale Superiore, Department of Mathematics, Piazza dei Cavalieri 7, Pisa, PI, Italy; e-mail: l.ambrosio@sns.it}, 
Gero Friesecke\footnote{TU M\"unchen, Zentrum Mathematik, Boltzmannstr.~3, D-85748 Garching, Germany; e-mail: gf@ma.tum.de}, 
Jannis Giannoulis\footnote{TU M\"unchen, Zentrum Mathematik, Boltzmannstr.~3, D-85748 Garching, Germany; e-mail: giannoulis@ma.tum.de}}
\\[5mm]
{July 7, 2009}
\end{center}
%
%
%

\begin{abstract}
We present a rigorous derivation of classical molecular dynamics (MD) from quantum
molecular dynamics (QMD) that applies to the 
standard Hamiltonians of molecular physics with Coulomb interactions. The derivation
is valid away from possible electronic eigenvalue crossings.
%
%
\end{abstract}

\section{Introduction}

A basic mathematical formulation of the passage from quantum to classical mechanics, following
the ideas of Eugene Wigner \cite{Wigner32}, is due to Lions and Paul \cite{LionsPaul93} and
G\'erard \cite{GerardWigner}: if $\{\Psi_\eps\}$ is a sequence of solutions to a semiclassically
scaled Schr\"odinger equation with smooth potential, then the sequence of corresponding Wigner transforms converges (up
to subsequences) to a solution of the Liouville equation, i.e. the transport equation of
the underlying classical dynamics. For a closely related mathematical approach to semiclassical
limits, which goes back to Egorov and is based on Weyl quantization and Moyal calculus, see e.g. \cite{Robert,
Martinez02}. 

Due to the reliance on smooth potentials, and in particular on the 
existence and uniqueness of trajectories of the classical dynamics, 
these results are not directly applicable when one tries to 
derive classical molecular dynamics
(MD) from Born-Oppenheimer quantum molecular dynamics (QMD). 
By the latter, one means 
quantum dynamics of the molecule's atomic nuclei in the exact non-relativistic Born-Oppenheimer
potential energy surface given by the ground state eigenvalue of the electronic Hamiltonian with
Coulomb interactions. The limit where
the natural small parameter in QMD, the ratio of electronic to nuclear mass $m_e/m_n=:\eps^2$, tends to zero,
has the structure of a semiclassical limit. (The physical value of this parameter is $\sim 1/2000$ for hydrogen, 
and even less for the other atoms.)
However, the potential energy surface of QMD
is not even continuous, because it always contains Coulomb singularities due to nuclei-nuclei repulsion;
in addition it can have cone-type singularities at electronic eigenvalue crossings. 

Here we present a rigorous derivation of MD from QMD in the limit of small mass ratio
that is applicable to the exact Born-Oppenheimer potential energy surface with Coulomb interactions. Our result
is valid away from eigenvalue crossings.
This is done by extending the approach 
of Lions and Paul \cite{LionsPaul93} to an appropriate class of non-smooth potentials. The main
technical novelty is a {  non-}concentration 
estimate on the set of Coulomb singularities, which allows
to show, in particular, that the singular term $\nabla U$ which appears in the Liouville equation
lies in $\L^1$ with respect to the limiting Wigner measure, thereby guaranteeing that the weak formulation
of the Liouville equation continues to make sense. 

Our methods do not seem to allow to analyse 
the limit dynamics at eigenvalue crossings,
since no analogon of our Coulombic non-concentration estimate is available. In fact, the physically
correct starting point to investigate what happens at crossings would not be QMD, as the Born-Oppenheimer
approximation underlying QMD also breaks down (see \cite{NeWi, Ze} for earliest insights, 
\cite{Hagedorn94} for a first rigorous account, and e.g. 
\cite{FeGe02, CdV03, Lasser04, LasserTeufel05} for recent results). 

In the remainder of this Introduction we describe our main result precisely. 
\\[2mm]
{\bf Quantum molecular dynamics}
To simplify matters we assume that all nuclei have equal mass.
In atomic units $(m_{e}=|e|=\hbar=1)$, non-relativistic Born-Oppenheimer
quantum molecular dynamics is given by the time-dependent Schr\"odinger equation
\begin{equation}\tag{SE} \label{SE}
  \begin{cases}
   i\eps\partial_t\Psieps(\cdot,t) =\Heps\Psieps(\cdot,t)
   \quad\text{for $t\in\R$},
 \\  
   \Psieps(\cdot,0)=\Psieps^0, 
  \end{cases}
\end{equation}
with Hamiltonian
\be \label{Heps} 
    \Heps = -\frac{\epsilon^2}2 \Delta + U,  
\ee
where $\Psieps(\cdot, t)\in \L^2(\R^d;\C)$ is the wavefunction of the nuclei 
at time $t$,
$d=3M$ ($M=$ number of nuclei), $\Delta$ is the Laplacian on $\R^d$, $\epsilon:=(m_e/m_n)^{1/2}$ is the (dimensionless) 
small parameter already discussed above (where we have assumed for simplicity that all nuclei have equal mass), and 
$U\, : \, \R^d\to\R$ is the Born-Oppenheimer
ground state potential energy surface obtained by minimization over electronic
states (see e.g.\ \cite{SO93}). The precise definition of $U$ is as follows.
Let $Z_1,..,Z_M\in\N$
and $R_1,..,R_M\in\R^3$ denote the charges and positions of the nuclei, and
let $N$ denote the number of electrons in the system (usually $N=\sum_{\alpha=1}^MZ_
\alpha$). Then for $x=(R_1,\ldots,R_M)\in\R^d$
\begin{align}
& U = E_{e\ell}+V_{nn}, 
\qquad E_{e\ell}(x) = \inf_{\psi}\langle\psi, H_{e\ell}(x)\psi\rangle, 
\label{pes}
\\& V_{nn}(x) = \sum_{1\le \alpha<\beta\le M} 
\frac{Z_\alpha Z_\beta}{|R_\alpha-R_\beta|}, 
\qquad H_{e\ell}(x) = \sum_{i=1}^N\Bigl(-\frac12\Delta_{r_i} 
- \sum_{\alpha=1}^M\frac{Z_\alpha}{|r_i-R_\alpha|}\Bigr) 
+ \sum_{1\leq i < j \leq N}\frac{1}{|r_i-r_j|}. 
\notag
\end{align}
Here the $r_i\in\R^3$ denote electronic coordinates and
the infimum is taken over the usual subset of $\L^2((\R^3\times\Z_2)^N;\C)$
of normalized, antisymmetric electronic states belonging to the domain
$\H^2((\R^3\times\Z_2)^N;\C)$ of $H_{e\ell}(x)$.
Physically, $E_{e\ell}$ is the electronic part of the energy, consisting of kinetic energy of the
electrons, electron-nuclei attraction, and electron repulsion; this part depends indirectly on
the positions $R_\alpha$ of the nuclei 
since these appear as parameters in the electronic Hamiltonian, 
and can be shown to be bounded and globally Lipschitz (cf.\ \cite{GeroQMBook}) 
although it is not elementary to see this. 
In case $N\le \sum_{\alpha=1}^M Z_\alpha$, Zhislin's theorem 
(see \cite{Fr03} for a short proof) 
says that the infimum in (\ref{pes}) is actually attained, the minimum value being an isolated eigenvalue of finite multiplicity of 
$H_{e\ell}$. $V_{nn}$ is the
direct electrostatic interaction energy between the nuclei, and is the origin of the discontinuities of $U$.

We remark that quantum molecular dynamics (\ref{SE}), \eqref{Heps}, 
(\ref{pes}), 
which is taken as starting point here,
itself already constitutes an approximation to full Schr\"odinger dynamics for electrons and nuclei.
Its rigorous justification constitutes an interesting problem in its own right; for a
comprehensive treatment in the case of smooth interactions and absence of electronic eigenvalue
crossings see \cite{TeufelHabil, PanatiSpohnTeufel}. 

The potential (\ref{pes}) satisfies the standard {  Kato-type} condition
\be \label{katocdn}
   U = U_b + U_s, \;\;\; U_b\in \L^\infty(\R^d), \;\;\;
   U_s(x) = \sum_{1\le\alpha<\beta\le M}V_{\alpha\beta}(R_\alpha-R_\beta), 
   \quad V_{\alpha\beta}\in \L^{2}(\R^3)+\L^\infty(\R^3).
\ee
For such potentials, the operator $\Heps$ is self-adjoint on $\L^2(\R^d)$ 
with domain ${\cal D}(\Heps)=\H^2(\R^d)$, cf.\ \cite{Kato}. 
By standard results on the unitary group generated by a self-adjoint operator, for any initial state 
$\Psieps^0\in {\cal D}(\Heps)$
this equation has a unique solution $\Psieps\in \Ce(\R;\H^2(\R^d))\cap \Ce^1(\R;\L^2(\R^d))$, the solution operator
$U_\eps(t) \, : \, \Psieps^0\mapsto\Psieps(\cdot,t)$ being unitary. 
In particular, 
\be \label{estSchroed} 
               \|\Psieps(\cdot,t)\|=\|\Psieps^0\| \mbox{ for all }t\in\R,
\ee
where here and below $\|\cdot \|$ denotes the $\L^2(\R^d)$ norm. 
\\[2mm]
{\bf The Wigner picture}
Given the state $\Psieps(\cdot,t)\in \L^2(\R^d;\C)$ of the system at time $t$,
define the associated Wigner function of lengthscale $\eps$ on $\R^d\times\R^d$, 
\begin{eqnarray}
  \Weps(x,p,t) & = & \frac{1}{(2\pi)^d}\int_{\R^d}
  \Psieps(x+\frac{\eps y}{2},t)\overline{\Psieps(x-\frac{\eps y}{2},t)} 
e^{-ip\cdot y}\d y \label{defW} \\
  & = & \frac{1}{(2\pi\eps)^d}\int_{\R^d} \Psieps(x+\frac{y}{2},t)
        \overline{\Psieps(x-\frac{y}{2},t)} e^{-ip\cdot {y}/{\eps}}\d y. 
\nonumber
\end{eqnarray}
Note that the integrand belongs to $\L^1(\R^d_y)$, so $\Weps$ is well defined
for a.e.\ $x\in\R^d$. In fact, the integrand is continuous in $x$ with respect to the  $\L^1(\R^d_y)$-norm, and hence $\Weps$ is continuous in $x$.
Roughly speaking, $\Weps$ is a joint position and 
momentum density of the system. Warning: $\Weps$ is not nonnegative 
except in the limit $\eps\to 0$, but at least its marginals are,
\begin{eqnarray}
  \int_{\R^d}\Weps(x,p,t) \d p & = & \Bigl|\Psieps(x,t)\Bigr|^2 \;\;(\mbox{position density}), \label{marginals}\\
  \int_{\R^d}\Weps(x,p,t) \d x & = & \Bigl|\frac{1}{(2\pi\eps)^{d/2}}
\underbrace{\int_{\R^d}e^{-ip\cdot{x}/{\eps}}\Psieps(x,t)\d x
}_{=:(\calF\Psieps)(p/\eps,t)}
\Bigr|^2 \;\;(\mbox{momentum density}).
\nonumber
\end{eqnarray}
Here and below, $\calF\phi$ denotes the (standard, not scaled) Fourier transform of the function $\phi$. 

When $\Psieps$ satisfies \eqref{SE}, its Wigner function satisfies
 \begin{align}\tag{WE} \label{WE}
  & \partial_t \Weps = - p\cdot \nabla_x \Weps + f_\eps,  \\
  & f_\eps(x,p,t) = -\frac{i}{(2\pi)^d} \int_{\R^d}
           \frac{U(x+\frac{\eps y}2) - U(x-\frac{\eps y}2)}{\eps} 
  \Psieps(x+\frac{\eps y}2,t)\overline{\Psieps(x-\frac{\eps y}2,t)} 
  e^{-ip\cdot y} \d y. \nonumber
 \end{align}
Formally, this follows from a lengthy but elementary calculation which goes back to Wigner,
see Section \ref{S:Wigner} below. In this section we also introduce a suitable function space setting in which 
the calculation becomes rigorous for general self-adjoint Hamiltonians of type (\ref{Heps}).   
We call eq. (\ref{WE}) the \emph{Wigner equation}. Note that it contains no modification
of quantum dynamics (\ref{SE}), but is just a different mathematical formulation of it.
\\[2mm]
{\bf Limit dynamics} 
>From now on we focus on the specific potential energy surface (\ref{pes}). 
In the limit $\eps\to 0$, the difference quotient in the potential term satisfies
$$
   \frac{U(x+\frac{\eps y}2) - U(x-\frac{\eps y}2)}{\eps} \to \nabla U(x)\cdot y \;\;\; \mbox{a.e.}
$$
(due to the Lipschitz continuity of $E_{e\ell}$ and the fact that the singular set of $V_{nn}$ is of
measure zero). To understand what happens with eq. (\ref{WE}) in the limit, it is useful to split $f_\eps$
into a term containing $\nabla U(x)\cdot y$ and a term containing the difference quotient of $U$ minus
its limit. The first term simplifies due to $y e^{-ip\cdot y}= i\nabla_p e^{-ip\cdot y}$, giving
$$
  f_\eps = \nabla U(x)\cdot \nabla_p\Weps + g_\eps,
$$
where 
$$
  g_\eps = -\frac{i}{(2\pi)^d} \int_{\R^d}
          \Bigl[ \frac{U(x+\frac{\eps y}2) - U(x-\frac{\eps y}2)}{\eps} - \nabla U(x)\cdot y\Bigr] 
  \Psieps(x+\frac{\eps y}2,t)\overline{\Psieps(x-\frac{\eps y}2,t)}
  e^{-ip\cdot y} \d y.
$$
Formally, passing to the limit in (\ref{WE}) and assuming that $g_\eps$ tends to zero,
we obtain the \emph{Liouville equation} 
\be
\partial_t W  =  - p\cdot \nabla_x W + \nabla U(x)\cdot \nabla_p W.
\tag{LE}\label{LE}
\ee
This is the transport equation for classical molecular dynamics in $\R^d\times\R^d$ with potential $U$,
\begin{equation}
\dot{x} =  p, \;\;\;
\dot{p} = -\nabla U(x). 
\tag{MD}\label{MD}
\end{equation}
Moreover, assuming the initial data to \eqref{SE} to be normalized, 
by (\ref{marginals}) 
one expects
$\int_{\R^{2d}} \dd W(t) = 1$, 
i.e. $W(t)$ should be a probability measure on phase space for all $t\in\R$. 

The main difficulties in making this rigorous for rough potentials
lie in (a) justifying the existence of a limiting probability measure on phase 
space for all times and general initial data
which are not restricted to `avoid' the singularities, 
(b) justifying that $g_\eps$ goes to zero. The latter issue arises because
when the potential $U$ is not everywhere differentiable, the term in square brackets does not go to zero
for every $x$, let alone locally uniformly. On the other hand, the remaining part of the integrand can
concentrate on individual positions $x$ when standard semiclassical wave packets such as
$$
    \Psieps^0(x) = \eps^{-\frac{\alpha d}2} e^{i\frac{p_0}{\eps}\cdot x} \phi(\mbox{$\frac{x-x_0}{\eps^\alpha}$}), \;\;\; 0<\alpha<1, \;\;\; \|\phi\|=1,
$$
are under consideration, whose Wigner function converges to $\delta_{(x_0,p_0)}$. Thus the only viable strategy
appears to be to establish that the Wigner function does not charge the singular set in the limit. 
\\[2mm]
{\bf Main result} 
Our rigorous result achieves goal (a) in the desired generality, and goal (b) away from possible crossings, establishing
in particular that the Liouville equation (\ref{LE}) remains valid across
Coulomb singularities.
{  In order to formulate our result we need the following definition.}
\begin{definition} 
A sequence $\{\mu_\eps\}$ of nonnegative Radon
measures on $\R^d$ is called {\em tight} if 
\begin{equation*}
   \lim_{R\to\infty} \limsup_{\eps\to0}\int_{|x|>R}\d \mu_\eps = 0. 
\end{equation*}
\end{definition}%
\begin{theorem}\label{limitequation} 
Suppose $U:\R^d\to\R$ is the Born-Oppenheimer potential energy surface (\ref{pes}) of any 
molecule, or more generally $U=U_b+U_s$
with $U_b\in \W^{1,\infty}(\R^d)$, 
\begin{equation}\label{Us}
U_s(x)=\sum_{1\le\alpha<\beta\le M}\frac{C_{\alpha\beta}}{|R_\alpha-R_\beta|},
\quad  C_{\alpha\beta}\ge0,  \quad x=(R_1,\ldots,R_M)\in\R^d.
\end{equation}
Let $\{\Psieps^0\}_{\eps>0}$ be a sequence of initial data such that
$\Psieps\in\H^2(\R^d)$, $\|\Psieps^0\|=1$,
$\|H_\eps\Psieps^0\|\le c$ for some constant $c$ independent of $\eps$, $\{|\Psieps^0|^2\}$ tight. 
Let $\Psieps\in \Ce(\R;\H^2(\R^d))\cap \Ce^1(\R;\L^2(\R^d))$ 
be the corresponding 
solutions to the time-dependent Schr\"odinger equation \eqref{SE}, 
and let $\Weps$ be their Wigner transforms \eqref{defW}.
Then:
\\[1mm]
(i) (Compactness) For a subsequence,  $\Weps\weakly W$ in $\mathcal{D}^\prime(\R^{2d+1})$. 
\\[1mm]
(ii) (Existence of a limiting probability measure on phase space) 
$W\in \Ce_{weak*}(\R;\calM(\R^{2d}))$, and $W(t)$ is a probability measure for all $t$, that is to say $W(t)\ge 0$ and
$\int_{\R^{2d}}\d W(t)=1$.
\\[1mm]
(iii) (No-concentration estimate at Coulomb singularities) 
For all $t$ we have $\nabla U_s\in\L^1(\d W (t))$, and 
$W(t)(\calS\times\R^d)=0$, where $\calS$ is the singular set
$\{x=(R_1,\ldots,R_M)\in\R^{3M}\ |\ R_\alpha=R_\beta$ for some  
$\alpha\neq\beta$ with $C_{\alpha\beta}\neq 0\}$.
\\[1mm] 
(iv) (Limit equation) If $\Omega\subseteq\R^d$ is any open set such that $U_b\in\Ce^1(\Omega)$,
 then  
$W$ is a global weak solution of the Liouville equation \eqref{LE} 
on $\Omega\times\R^d\times\R$, that is to say
 \be \label{weakLE}
    \int_\R\int_{\R^{2d}}\Bigl(\partial_t + p\cdot \nabla_x - 
 \nabla U(x)\cdot \nabla_p\Bigr)\phi(x,p,t) \d W(t)\d t 
    =0
 \ee
 for all $\phi\in \Ce_0^\infty(\Omega\times\R^d\times\R)$. 
\end{theorem}
{\bf Remarks}
1) The regularity requirement $U_b\in \Ce^1$ in (iv) is minimal in order for the weak Liouville equation (\ref{weakLE}) 
to make sense for general measure valued
solutions $W$. This is due to the appearance of $\nabla U$ inside the integral with respect to the measure $\d W(t)$.
Note however that in case of (\ref{pes}) this narrowly excludes eigenvalue
crossings, as seen from the 2D matrix example $$H_{e\ell}=\begin{pmatrix}\rho_1 & -\rho_2 \\ \rho_2 & \rho_1\end{pmatrix},$$ 
whose ground state eigenvalue equals $-|\rho|$, and is hence Lipschitz but not $\Ce^1$. For interesting model problems with scalar or vector-valued potentials
in which the behaviour of Wigner measures past discontinuities of $\nabla U$ can be analysed for suitable classes
of initial data see \cite{Keraani05} and \cite{FeGe02, Lasser04, LasserTeufel05}.  
\\[2mm] 
2) The assumptions on the potential $U$ are far weaker than those needed for uniqueness
of the Hamiltonian ODE \eqref{MD} underlying the limit equation. Recall that the standard condition guaranteeing
uniqueness for ODE's $\dot{z}=f(z)$ is boundedness of the gradient of the vector field $f$, which in the
case of \eqref{MD}  means boundedness of the {\it second}, not the first gradient of $U$.
Our assumptions on $U$ are also weaker than those under which uniqueness
for weak ($\L^p$) solutions to \eqref{LE} is known. 
The recent nontrivial uniqueness results for transport
equations (\cite{Ambrosio04}) require $f\in \mathrm{BV}$, i.e. in case of \eqref{MD}, $\nabla U\in \mathrm{BV}$ 
(for recent refinements see \cite{BC09}, \cite{AGS08}). 
Interestingly, however, the latter
requirement, while violated by the Coulombic part $U_s=V_{nn}$ in (\ref{pes}), would be met by the model eigenvalue crossing in Remark 1),
and expected to be met by the electronic part $U_b=E_{e\ell}$ in (\ref{pes}). We plan to address uniqueness in a separate publication \cite{AFFG09}.  
\\[2mm]
3) The higher integrability result in (iii) that $\nabla U_s\in \L^1(\d W(t))$, which is essential for making sense of the limit equation at Coulomb singularities,
requires a quantitative no-concentration estimate of form 
\be \label{noconcest}
             \int_{|R_\alpha-R_\beta|<\delta} \d W(t) = O(\delta^2) \mbox{ as }\delta\to 0, 
\ee
for any $\alpha\neq\beta$ with $C_{\alpha\beta}\neq 0$.
This is because $|\nabla_{R_\alpha}(1/|R_\alpha-R_\beta|)|=1/|R_{\alpha}-R_\beta|^2\ge 1/\delta^2$ in $|R_\alpha-R_\beta|<\delta$. 
We do not think that the validity of such an
estimate is obvious, the naively expected bound only being $O(\delta)$ instead of $O(\delta^2)$ 
(on grounds of the potential energy term $\int U_s |\Psieps|^2$,
which can be controlled independently of $\eps$ and $t$ by energy conservation, only containing the weaker singularities $1/|R_{\alpha}-R_\beta|$).
See Sections 4 and 5 for the proof of (\ref{noconcest}). 
\\[2mm] 
4) In the special case $U_s=0$, $U_b\in \W^{2,\infty}(\R^d)$, (iv) holds with $\Omega=\R^d$,
so Theorem \ref{limitequation} recovers
the result of Lions and Paul \cite[Th\'eor\`eme {IV}.1.1)]{LionsPaul93}.
\\[2mm]
5) In the special case of the potential energy surface (\ref{pes}) of the H$_2$ molecule 
($M=2$, $N=2$, $Z_1=Z_2=1$), it is known that the ground state eigenvalue of the electronic Hamiltonian
is nondegenerate. It then follows from a result of Hunziker \cite{Hun86} that $U_b$ is analytic, and in particular
$\Ce^1$, on $\Omega=\R^6\backslash{\cal S}$, 
and Theorem \ref{limitequation} justifies classical molecular dynamics globally. 
\\[2mm]
6) Also, more can be said about the set $\Omega$ in (iv) 
in the case of the potential energy surface (\ref{pes}) of a general neutral or
positively charged dimer ($M=2$, $N\le Z_1+Z_2$).
By invariance of the electronic Hamiltonian
$H_{e\ell}$ and the nuclei-nuclei interaction $V_{nn}$ in (\ref{pes}) under simultaneous rotation and translation of
all particles, we have that $U(R_1,R_2)=u(|R_1-R_2|)$, that is to say the potential is a function of a single parameter,
internuclear distance. It then follows by combining the result of Hunziker \cite{Hun86}
with a classical result of Kato on analyticity of eigenvalues of analytic one-parameter families of
Hamiltonians \cite{Kat95} that we may take $\Omega=\R^6\backslash({\cal S}\cup{\cal C})$, where ${\cal S}=\{R_1=R_2\}$ 
is the set of Coulomb singularities introduced in part (iii) of the theorem, and ${\cal
C}=\bigcup_j\{|R_1-R_2|=c_j\}$,
the $c_j$ being the (possibly empty) discrete subset of $\R^+$ of interatomic distances at which
the lowest two eigenvalues of the electronic Hamiltonian cross.
Note in particular that ${\cal S}\cup{\cal C}$ is a closed set of measure zero; hence Theorem \ref{limitequation} justifies
the Liouville equation on an open set of full measure. We expect that when $U$ is given by (\ref{pes}), it is always smooth
on an open set of full measure. Note however that such a result does not follow solely from consideration
of non-degenerate eigenvalues as in \cite{Hun86}.
\\[2mm]
7) As shown below (Lemma \ref{L:timedependent}), the weak convergence in (i) also holds in the stronger spaces 
$\L^\infty_{weak\ast}(\R;\calA^\prime)$ 
and $\Ce_{weak\ast,loc}(\R;\calA^\prime)$,
with $\calA$ being the Banach space defined in \eqref{star}.
\section{Wigner-transformed quantum dynamics}
\label{S:Wigner}
We now make precise in an appropriate function space setting
the well known fact (discussed informally in the Introduction)
that the Wigner transform takes
solutions of the Schr\"odinger equation \eqref{SE} 
to solutions of the Wigner equation (\ref{WE}). 
Our choice of spaces is convenient
for our goal to study the limit dynamics for rough potentials. 
In other contexts other function spaces have been considered 
\cite{Markowich}.

We begin with the well known formal derivation (assuming 
that the wavefunction is smooth and rapidly decaying
for all $t$).
\\[2mm]
{\bf Formal derivation} Let $\Psi_\eps$ be a solution to (\ref{SE}), 
and let $W_\eps$ denote its Wigner transform (\ref{defW}). 
Differentiating the latter 
with respect to $t$, 
we obtain
\be \label{twoa}
  \partial_t\Weps(x,p,t)  =  \frac{1}{(2\pi)^d}\int_{\R^d}
  \Bigl[\bigl(\partial_t\Psieps(x+\frac{\eps y}{2},t)\bigr)
\overline{\Psieps(x-\frac{\eps y}{2},t)} 
     + \Psieps(x+\frac{\eps y}{2},t)\overline{\partial_t\Psieps(x-\frac{\eps y}{2},t)}
     \Bigr] e^{-ip\cdot y}\d y. 
\ee
{  
By (\ref{SE}) and $\Delta_{\pm\eps y/2} = (4/{\eps}^2)\Delta_y$
this is equivalent to
\begin{multline}
  \partial_t\Weps(x,p,t)   = f_\eps(x,p,t)
\\ +
\frac{2i}{\eps (2\pi)^d} 
\int_{\R^d}
  \Bigl[\bigl(\Delta_y \Psieps(x+\frac{\eps y}2,t)\bigr)
\overline{\Psieps(x-\frac{\eps y}2,t)} 
-\Psieps(x+\frac{\eps y}2,t) \overline{\Delta_y\Psieps(x-\frac{\eps y}2,t)}
  \Bigr] e^{-ip\cdot y}\d y
\label{evolWeps1}
\end{multline}
with $f_\eps$ as in \eqref{WE}.
The Laplacian terms can be simplified via the formula 
$(\Delta a)\overline{b} - a\overline{\Delta b} 
=\div(\nabla a\cdot \overline{b} 
- a \cdot\nabla \overline{b})$, an integration by parts, and the formula
$\nabla_y=\pm (\eps/2)\nabla_{\pm \eps y/2}$, 
whence the integral in \eqref{evolWeps1} becomes
\begin{align}
& \int_{\R^d}
    \div\Bigl[\bigl(\nabla_y \Psieps(x+\frac{\eps y}2,t)\bigr)
\overline{\Psieps(x-\frac{\eps y}2,t)} 
    - \Psieps(x+\frac{\eps y}2) \overline{\nabla_y\Psieps(x-\frac{\eps y}2,t)}
\Bigr]
  e^{-ip\cdot y}\d y \label{evolWeps2}\\
\displaybreak[0]  &  = i p\cdot \int_{\R^d}
    \Bigl[ \bigl(\nabla_y \Psieps(x+\frac{\eps y}2,t)\bigr) 
\overline{\Psieps(x-\frac{\eps y}2,t)}
- \Psieps(x+\frac{\eps y}2,t)  \overline{\nabla_y\Psieps(x-\frac{\eps y}2,t)}
        \Bigr]
e^{-ip\cdot y}
\d y \label{evolWeps3} \\
\displaybreak[0]   &  = 
\frac{\eps i}{2}\,p\cdot 
\int_{\R^d} 
\nabla_x\Bigl[\Psieps(x+\frac{\eps y}2,t)\overline{\Psieps(x-\frac{\eps y}2,t)}
\Bigr]
e^{-ip\cdot y} 
\d y
    \label{evolWeps4}\\
\displaybreak[0] & = 
\frac{\eps(2\pi)^d}{2i}
\,\bigl(- p\cdot \nabla_x \Weps(x,p,t)\bigr).
\notag
\end{align}
Substituting this expression into 
(\ref{evolWeps1}), 
we obtain (\ref{WE}).
}
\\[2mm]
{\bf Rigorous derivation} 
In the sequel, position coordinates in $\R^d$, $d=3M$, are denoted by
$x=(x_1,..,x_d)=(R_1,..,R_M)\in\R^d$,
$x_i\in\R$, $R_\alpha\in\R^3$.
\begin{lemma}\label{L:rigorousWigner}
The Wigner transform (\ref{defW}) of any solution 
$\Psieps\in \Ce(\R;\H^2(\R^d))\cap \Ce^1(\R;\L^2(\R^d))$
of the Schr\"odinger equation \eqref{SE} 
with $U$ as in \eqref{katocdn} 
satisfies 
\begin{equation}\label{wignerregularity}
   \Weps\in \Ce^1(\R;\L^\infty (\R^{2d})), \;
   \frac{\partial}{\partial x_i}\Weps,
\ \frac{\partial^2}{\partial x_i\partial x_j}\Weps, 
\mbox{ and }f_\eps\in \Ce(\R;\L^\infty (\R^{2d}))
\quad\text{for all}\quad i,j=1,\ldots,d
\end{equation}
and solves the Wigner equation (\ref{WE}).
\end{lemma}
To obtain an effortless proof, the idea is to express all terms under investigation 
with the help of the following bilinear map which extends the quadratic map 
{  $\Psieps\mapsto\Weps$ introduced in \eqref{defW}:} 
$$
   \Feps(\Psi,\chi) := \frac{1}{(2\pi)^d}\int_{\R^d} \Psi(x+\mbox{$\frac{\eps y}{2}$}) \overline{\chi(x-\mbox{$\frac{\eps y}{2}$})}
   e^{-ip\cdot y}\d y.
$$
\begin{lemma} \label{L:bilinear} 
The map $(\Psi,\chi)\mapsto \Feps(\Psi,\chi)$ is a continuous map from 
$\L^2(\R^d)\times \L^2(\R^d)$ to $\L^\infty(\R^{2d})$. 
In particular, the map $\Psieps\mapsto\Weps=\Feps(\Psieps,\Psieps)$ 
is a continuous map from $\L^2(\R^d)$ to $\L^\infty(\R^{2d})$. 
\end{lemma}
{\bf Proof of Lemma \ref{L:bilinear}} Let $\Psi$, $\Psi'$, $\chi$, $\chi'\in \L^2(\R^d)$, $W=\Feps(\Psi,\chi)$, 
$W'=\Feps(\Psi',\chi')$. Then 
\begin{align*}
  &  |W(x,p)-W'(x,p)| \\
  & =  \frac1{(2\pi)^d}\left| \int_{\R^d} \Bigl( (\Psi-\Psi')(x+\mbox{$\frac{\eps y}{2}$})\overline{\chi(x-\mbox{$\frac{\eps y}{2}$})} 
                            + \Psi'(x+\mbox{$\frac{\eps y}{2}$})
  \overline{(\chi-\chi')(x-\mbox{$\frac{\eps y}{2}$})} \Bigr) e^{-ip\cdot y}
  \d y \right| \\
  &  \le  \frac1{(2\pi)^d}\Bigl( \|(\Psi - \Psi')(x+\mbox{$\frac{\eps \cdot }{2}$})\|\,
  \|\chi(x-\mbox{$\frac{\eps \cdot }{2}$})\|
        + \|\Psi'(x + \mbox{$\frac{\eps \cdot }{2}$})\|\,
\|(\chi-\chi')(x-\mbox{$\frac{\eps \cdot }{2}$})\|\Bigr) \\
  & =  \Bigl(\frac1{\eps\pi}\Bigr)^d\Bigl(\|\Psi-\Psi'\| \, \|\chi\| + \|\Psi'\| \, \|\chi-\chi'\|\Bigr)
\end{align*}
for all $x$ and $p$. Taking $\Psi'=\chi'=0$ shows $W\in \L^\infty(\R^{2d})$, 
and the estimate above establishes the 
asserted continuity of $\Feps$.\hfill$\square$\\[2mm]
{\bf Proof of Lemma \ref{L:rigorousWigner}} 
First, we claim that $\Weps\in\Ce(\R;\L^\infty(\R^d))$. 
This is immediate from $\Psieps\in\Ce(\R;\L^2(\R^d))$ and Lemma \ref{L:bilinear}.

Next, we investigate the terms $\frac{\partial}{\partial x_i}\Weps$, $\frac{\partial^2}{\partial x_i\partial x_j}\Weps$, and $f_\eps$. 
The underlying terms $\frac{\partial}{\partial x_i}\Psieps$,
$\frac{\partial^2}{\partial x_i\partial x_j}\Psieps$ 
and $U\Psieps$ are in $\Ce(\R;\L^2(\R^d))$,
because 
$\Psieps\in \Ce(\R;\H^2(\R^d))$ and the operators 
$\frac{\partial}{\partial x_i}$, $\frac{\partial^2}{\partial x_i\partial x_j}$ 
and $U=\Heps - \frac{\eps^2}{2}\Delta$ are continuous
maps from $\H^2(\R^d)$ to $\L^2(\R^d)$. 
Note now that we have the following representations with the help of the 
bilinear map $\Feps$:
\begin{eqnarray*}
   & & \frac{\partial}{\partial x_i}\Weps =\textstyle  \Feps\Bigl(\frac{\partial}{\partial x_i}\Psieps,\Psieps\Bigr) + \Feps\Bigl(\Psieps, \frac{\partial}{\partial x_i}\Psieps\Bigr), \\
   & & \frac{\partial^2}{\partial x_i\partial x_j}\Weps 
\textstyle = \Feps\Bigl(\frac{\partial^2}{\partial x_i\partial x_j}\Psieps,\Psieps\Bigr) 
+ \Feps\Bigl(\frac{\partial}{\partial x_i}\Psieps,
\frac{\partial}{\partial x_j}\Psieps\Bigr) 
+ \Feps\Bigl(\frac{\partial}{\partial x_j}\Psieps,
\frac{\partial}{\partial x_i}\Psieps\Bigr) 
+ \Feps\Bigl(\Psieps,\frac{\partial^2}{\partial x_i\partial x_j}\Psieps\Bigr), 
\\
   & & f_\eps = -\frac{i}{\eps}\Bigl(\Feps(U\Psieps,\Psieps) - \Feps(\Psieps,U\Psieps)\Bigr).
\end{eqnarray*}
It now follows from Lemma \ref{L:bilinear} that these terms are in 
$\Ce(\R;\L^\infty(\R^{2d}))$. 

It remains to show that $\Weps\in \Ce^1(\R;\L^\infty(\R^{2d}))$, 
with derivative 
$\partial_t\Weps$ 
given by eq. (\ref{WE}). 
First we show continuous differentiability with respect to time. We have
\begin{eqnarray*}
  \frac{\Weps(\cdot,t+h)-\Weps(\cdot,t)}{h} &=& \frac{\Feps(\Psieps(\cdot,t+h),\Psieps(\cdot,t+h)) - \Feps(\Psieps(\cdot,t),\Psieps(\cdot,t))}{h} \\
    &=& \Feps\Bigl(\mbox{$\frac{\Psieps(\cdot,t+h)-\Psieps(\cdot,t)}{h}$}, \Psieps(\cdot,t+h)\Bigr) +
        \Feps\Bigl(\Psieps(\cdot,t), \mbox{$\frac{\Psieps(\cdot,t+h)-\Psieps(\cdot,t)}{h}$}\Bigr).
\end{eqnarray*}
Hence by the two convergences 
$\frac{\Psieps(\cdot,t+h)-\Psieps(\cdot,t)}{h}\to\partial_t\Psieps(\cdot,t)$ in $\L^2(\R^d)$ 
(from $\Psieps\in\Ce^1(\R;\L^2(\R^d))$)
and $\Psieps(\cdot,t+h)\to\Psieps(\cdot,t)$ in $\L^2(\R^d)$ 
(from $\Psieps\in\Ce(\R;\L^2(\R^d))$) 
and the continuity of $\Feps$, 
$$
   \frac{\Weps(\cdot,t+h)-\Weps(\cdot,t)}{h} \to \Feps(\partial_t\Psieps(\cdot,t),\Psieps(\cdot,t)) + \Feps(\Psieps(\cdot,t),\partial_t\Psieps(\cdot,t)) \mbox{ in }\L^\infty(\R^{2d})
   \mbox{ as }h\to 0.
$$
Consequently $t\mapsto\Weps(\cdot,t)$ is a differentiable map from $\R$ to $\L^\infty(\R^{2d})$, with derivative given by eq. (\ref{twoa}). 
Continuity in time of the derivative $\partial_t\Weps$, i.e.
the fact that $\Weps\in\Ce^1(\R;\L^\infty(\R^{2d}))$,  
now follows from (\ref{twoa}), $\Psieps\in\Ce(\R;\L^2(\R^d))$, $\partial_t\Psieps\in\Ce(\R;\L^2(\R^d))$, and -- one more time -- the continuity of $\Feps$
(Lemma \ref{L:bilinear}). 

We conclude the proof by showing that $W_\eps$ satisfies \eqref{WE}.
The formal derivation of this equation from eq. (\ref{twoa}) (which we have already established above)
has been performed by the calculations \eqref{evolWeps1}---\eqref{evolWeps4}. 
Here, we need only to justify these calculations rigorously. 
Eq.\ \eqref{evolWeps1} follows immediately from (\ref{twoa}) and \eqref{SE}.
Note that 
{  all four summands of the integrands on the RHS of eq.\ \eqref{evolWeps1}
(cf.\ also the definition of $f_\eps$ in \eqref{WE}),}
separately belong to
$\L^1(\R^d_y)$, for any $x$, $p$ and $t$, because 
$U(x\pm\frac{\eps y}2)\Psieps(x\pm\frac{\eps y}2,t)$ and 
$\Delta_x\Psieps(x\pm\frac{\eps y}2,t)$ 
belong to $\L^2(\R^d_y)$.
%
Eq.\ \eqref{evolWeps2} follows from the product rule for the Laplacian, and
eq.\ \eqref{evolWeps3} from the fact that the vector field inside the 
square brackets of \eqref{evolWeps2} belongs to the Sobolev space 
$\W^{1,1}(\R^d_y)=\W_0^{1,1}(\R^d_y)$ and 
$e^{-ip\cdot y}\in\W^{1,\infty}(\R^d_y)$, and that 
$\int_{\R^d}\div(v)\phi=-\int_{\R^d}v\cdot\nabla\phi$ for all 
$v\in\W_0^{1,1}(\R^d)$, $\phi\in\W^{1,\infty}(\R^d)$.
Finally,  eq.\ \eqref{evolWeps4} follows from an elementary change of
 variables,
concluding the proof. \hfill$\square$\\
%
%
%
%
%
\section{Time-dependent Wigner measures}

We first give a modest technical extension of the construction of
Wigner measures in \cite{LionsPaul93}. Instead of considering a sequence of wavefunctions at a fixed time $t$, 
we consider a sequence of continuous paths $\Psieps\in \Ce(\R;\L^2(\R^d))$ -- not required to satisfy any equation --
and show that under mild conditions these give rise to a continuous path 
$W\in \Ce_{weak*}(\R;{\calM}(\R^{2d}))$ of Wigner measures. See Lemma \ref{L:timedependent}. We then combine the lemma with careful a priori
estimates on the singular contributions to $f_\eps$ in eq. (\ref{WE}) and a lemma on propagation of tightness
under (\ref{SE})
to prove Theorem \ref{limitequation} (i) and (ii). An interesting feature of these proofs is
that, in contrast to the existing literature, they are extracted directly from the Schr\"odinger dynamics, 
without relying on a representation of the limit measure as push-forward of its initial data
under (\ref{MD}), which is not available here. 
As regards our proof in (ii) that $\int_{\R^{2d}} \d W(t) = 1$ for all $t$, our argument
via Schr\"odinger dynamics is inspired by Corollary 1 in \cite{FL03} (see also Proposition 4 in \cite{Lasser04}).

We begin by recalling the notion of weak convergence of Wigner transforms at
a fixed time $t$ to Wigner measures introduced in \cite{LionsPaul93}.
For an alternative construction of
Wigner measures see \cite{GerardWigner}. 
Closely related constructions 
are the microlocal defect measures of G\'erard \cite{GerardHmeasures},
motivated by questions in microlocal analysis, and the $H$-measures of Tartar 
\cite{Tartar90}, motivated by questions in
homogenization theory. 

Let $\calA$ denote the following Banach space
\begin{equation}\label{star}
   \calA:=\{\phi\in \Ce_0(\R^{2d}) \, | \, \|\phi\|_{\calA} := 
   \int_{\R^d}\sup_{x\in\R^d}|(\calF_p\phi)(x,y)|\d y <\infty \}.
\end{equation}
Here $\Ce_0(\R^{2d})$ is the usual space of continuous functions on $\R^{2d}$ 
tending to zero at infinity, and
$\calF_p\phi$ is the partial Fourier transform
$
   (\calF_p\phi)(x,y)=\int_{\R^d} e^{-ip\cdot y} \phi(x,p)  \d p.
$
Since $\calA$ is a dense subset of $\Ce_0(\R^{2d})$, its dual $\calA^\prime$ 
contains $\Ce_0^\prime(\R^{2d})={\cal M}(\R^{2d})$, 
the space of not necessarily nonnegative Radon measures on $\R^{2d}$ of finite mass. In particular, the
delta function $\delta_{(x_0,p_0)}$ centered at a single point $(x_0,p_0)$ in classical phase space belongs
to $\calA'$, and weak* convergence in $\calA'$ allows convergence of smeared-out Wigner functions coming
from a quantum state to a delta function on phase space (i.e., a ``classical'' state). 

The following basic facts were established by Lions and Paul.
\begin{lemma} \label{L:stationary} (Wigner measures) \cite{LionsPaul93} \\[1mm]
(i) (compactness) Let $\{\Psieps\}$ be any sequence in $\L^2(\R^d)$ such that
\be \label{psibound}
     \|\Psieps\|^2\le C
\ee
for some constant $C$ independent of $\epsilon$. Then the sequence of corresponding Wigner transforms $\{\Weps\}$
contains a subsequence $\{W_{\eps'}\}$ 
converging weak* in $\calA'$ to some $W\in\calA'$.  
\\[1mm]
(ii) (positivity) Any such limit $W\in\calA'$ satisfies $W\in\calM(\R^{2d})$, $W\ge 0$. In other words, $W$ is a nonnegative
Radon measure of finite mass.
\\[1mm]
(iii) (upper bound) Let $\{\Psi_{\eps'}\}$ 
be a further subsequence such that $|\Psi_{\epsilon'}|^2\weakstar\mu$ in 
$\calM(\R^{d})$. Then 
$$
  \int_{p\in\R^d} W(\cdot,\dd p) \le \mu.
$$
In particular, any limit $W$ as in (i) satisfies $\int_{\R^{2d}} \d W \le C$, 
with $C$ as in (\ref{psibound}). 
\\[1mm]
(iv) (preservation of mass) If $\|\Psieps\|^2=C$ for all $\eps$, 
and the sequences of position and momentum densities are both tight, that is to say 
\begin{equation*}
\limsup_{\eps\to0}\int_{|x|>R}\left|\Psieps(x)\right|^2\d x\to 0,\quad
\limsup_{\eps\to0}\int_{|p|>R}
\left|\frac{1}{(2\pi\eps)^{d/2}}(\calF\Psieps)\left(\frac{p}{\eps}\right)\right|^2\d p\to 0
\quad (R\to\infty),
\end{equation*}
then $\displaystyle\int_{\R^{2d}}\d W=C$.
In particular, in the case $C=1$, $W$ is a probability measure.
\end{lemma}
For future purposes we note 
that (i) is immediate from the Banach-Alaoglu theorem and the 
elementary estimate
\be
  \Bigl|\int_{\R^{2d}} \Weps\phi \, \d (x,p) \Bigr| 
  \le \frac1{(2\pi)^d} \|\Psieps\|^2 \|\phi\|_{\calA}
\quad\text{for $\phi\in\calA$}, 
\label{elest}
\ee
which implies 
\be \label{elesttwo}
  \|\Weps\|_{\calA'} = \sup_{\phi\in\calA\backslash\{0\} }\frac{\int \Weps \, \phi}{\|\phi\|_{\calA}} 
  \le \frac1{(2\pi)^d}\|\Psieps\|^2,
\ee
that is to say $\{\Weps\}$ is a bounded sequence in $\calA^\prime$. 
The proofs of (ii) and (iii) are less elementary, and require use of the Husimi transform; see \cite{LionsPaul93}. 
\\[2mm]
An analogue yielding time-continuous paths of Wigner measures is the following.
\begin{lemma}\label{L:timedependent} (Time-dependent Wigner measures) \\[1mm]
(i) (compactness) Let  $\{\Psieps\}$ be a sequence in $\Ce(\R;\L^2(\R^d))$ such that
\be \label{H1}
   \sup_{t\in\R}\|\Psieps(\cdot,t)\|^2\le C
\ee
for some constant independent of $\epsilon$, and let $\{\Weps\}$ 
be the sequence of associated Wigner functions. Then for
a subsequence, $\Weps\weakstar W$ weak* in $L^\infty(\R;\calA')$. 
\\[1mm]
(ii) Suppose in addition that for any test function $\phi\in \Ce_0^\infty(\R^{2d})$ the functions 
$$
  f_{\eps,\phi}(t) := \int_{\R^{2d}}\Weps(x,p,t)\, \phi(x,p)\d(x,p)
$$
are differentiable and satisfy
\be \label{H2}
    \sup_{t\in\R} |\mbox{$\frac{\dd}{\dd t}$}f_{\eps,\phi}(t)|\le C_\phi
\ee 
for some constant $C_\phi$ independent of $\epsilon$. Then, 
$W\in \Ce_{weak*}(\R;\calM(\R^{2d}))$, $W(t)\ge 0$ for all $t$,
and $\Weps(\cdot,t)\weakstar W(t)$ in $\calA^\prime$ for all $t$.
Moreover, the latter convergence is uniform on compact time intervals, 
i.e.,  for any test function 
$\phi\in\calA$ and any compact $I\subset\R$, $f_{\eps,\phi}(t)$
converges uniformly 
with respect to $t\in I$ to $\int_{\R^{2d}}\phi\,\d W(t)$.
\end{lemma}
%
%
%
Here (i) is a straightforward adaptation of the time-independent theory in \cite{LionsPaul93}.
The key point is the assertion in (ii) that the limit measure has slightly higher regularity in time than
na\"ively expected (continuous instead of $\L^\infty$). This allows, in particular, 
to make sense of initial values.        
\\[2mm]
{\bf Proof} 
The first part is an easy consequence of Lemma \ref{L:stationary},
cf.\ \eqref{elesttwo},  and assumption (\ref{H1}), which 
imply that $\{\Weps\}$ is bounded in $\L^\infty(\R;\calA')$. Since the latter is the dual of the separable space
$\L^1(\R;\calA)$, the assertion follows from the Banach-Alaoglu theorem. 

The second part is less trivial, and requires various approximation arguments. First, we test the weak* convergence of
$\{\Weps\}$ from (i) against tensor products $\phi(x,p)\chi(t)$ with $\phi\in\calA$, $\chi\in \L^1(\R)$. 
This gives
$$
   \int_\R f_{\epsilon,\phi}(t)\chi(t)\d t 
= \int_{\R^{2d+1}} W_\eps \, \phi\otimes\chi \d (x,p,t) \to
   \int_{\R}\int_{\R^{2d}}\phi\otimes\chi \d W(t) \d t 
= \int_{\R} f_\phi (t) \chi(t) \d t,
$$
where $f_\phi(t):=\int_{\R^{2d}}\phi\d W(t)$. 
Consequently $f_{\eps,\phi}\weakstar f_\phi$ in 
$\L^\infty(\R)$. 
%
%

Now let $\phi\in \Ce_0^\infty(\R^{2d})$. Then by assumption (\ref{H2}) and the
compact embedding 
{
$\W^{1,\infty}([-T,T])\hookrightarrow \Ce([-T,T])$, 
$f_{\eps,\phi}$ converges uniformly to $f_\phi$ 
on any compact interval $[-T,T]$. 
}
In particular $f_\phi$ is, as a uniform
limit of continuous functions, continuous, and $f_{\eps,\phi}(t)\to f_{\phi}(t)$ pointwise for all $t\in\R$, that is to say
\be \label{distconv}
   \Weps(\cdot, t) \rightharpoonup W(t) \mbox{ in }{\cal D}'(\R^{2d}) 
\mbox{ pointwise for all }t\in\R.
\ee
Now fix $t$. For a further subsequence which may depend on $t$, by Lemma \ref{L:stationary} 
$W_{\eps'}(\cdot,t)\weakstar \widetilde{W}(t)$
in $\calA'\subset {\cal D}'$. Together with (\ref{distconv}) this yields 
$W(t)=\widetilde{W}(t)$, as well as the
convergence $W_\eps(\cdot,t)\weakstar W(t)$ in ${\calA}'$ 
for the whole sequence. By Lemma \ref{L:stationary},
all remaining statements about 
$W(t)$ follow, except its asserted continuity in $t$. 

It remains to show the latter, i.e.
\be \label{Wcont}
   \int_{\R^{2d}} \phi \d W(t+h) \to \int_{\R^{2d}} \phi \d W(t) 
\;\;\; (h\to 0)\;\;\; \mbox{ for all }\phi\in \Ce_0(\R^{2d}). 
\ee
Given $\phi\in \Ce_0(\R^{2d})$ and $\delta>0$, 
by the density of $\Ce_0^\infty(\R^{2d})$ in $\Ce_0(\R^{2d})$ 
there exists $\phi_\delta\in \Ce_0^\infty(\R^{2d})$
such that $\|\phi_\delta - \phi\|_\infty<\delta$, 
where $\|\cdot\|_\infty$ denotes the norm of $\L^\infty(\R^{n})$. 
Consequently
\begin{eqnarray*}
  & & \Bigl| \int_{\R^{2d}} \phi \d W(t+h) - \int_{\R^{2d}}\phi\d W(t)\Bigr| 
\\
  & & \le \Bigl|\int_{\R^{2d}} \phi_\delta \d W(t+h)-\int_{\R^{2d}} \phi_\delta \d W(t)\Bigr| 
+ \Bigl|\int_{\R^{2d}}(\phi-\phi_\delta)\d W(t+h) \Bigr| 
         +  \Bigl|\int_{\R^{2d}}(\phi-\phi_\delta)\d W(t) \Bigr| \\
  & & \le   \Bigl|\int_{\R^{2d}} \phi_\delta \d W(t+h)-\int_{\R^{2d}}\phi_\delta\d W(t)\Bigr| 
+ \delta\Bigl(  \int_{\R^{2d}} \d W(t+h) +\int_{\R^{2d}}\d W(t)\Bigr).
\end{eqnarray*}
As $h\to 0$, the first term vanishes by the continuity of
$\int_{\R^{2d}}\phi_\delta \d W(t)=f_{\phi_\delta}(t)$ in $t$ 
(which was already established above). The second term stays bounded by
$2C\delta$ by {  Lemma \ref{L:stationary},} (iii), giving
$$
   \limsup_{h\to 0} \Bigl|\int_{\R^{2d}}\phi \d W(t+h) 
- \int_{\R^{2d}} \phi \d W(t)\Bigr| \le 2 C \delta.
$$
Since $\delta$ was arbitrary, the continuity assertion (\ref{Wcont}) follows,
and the proof of Lemma \ref{L:timedependent} is complete.\hfill$\square$\\[2mm]
{\bf Remark} 
Functional analytically, Lemma \ref{L:timedependent} says that 
under the assumptions (\ref{H1}), (\ref{H2}), $\{\Weps\}$ is
relatively compact in $\Ce_{weak*,\ell oc}(\R;\calA')$. 
This may be viewed as a weak-convergence variant of the
well known compactness lemma of J.~L.~Lions 
\cite[Chap.~1, Th\'eor\`eme 5.1]{JLL}, in which the
condition of boundedness of time derivatives in some Banach space has been replaced by condition (\ref{H2}) which
is related to a weak topology.
\\[5mm]
We close this section by applying the above lemma to prove the first two statements of Theorem \ref{limitequation}. 
\\[2mm]
{\bf Proof of Theorem \ref{limitequation} (i)} 
This follows from $\|\Psieps(t)\|=\|\Psieps^0\|=1$, 
Lemma \ref{L:timedependent} (i), 
and the fact that weak$^\ast$ convergence in $\L^\infty(\R;\calA^\prime)$ 
implies convergence in $\calD^\prime(\R^{2d+1})$.\hfill$\square$ 
\\[2mm]
{\bf Proof of Theorem \ref{limitequation} (ii)} 
First we prove that $W\in \Ce_{weak*}(\R;\calM(\R^{2d}))$, $W(t)\ge 0$.
By the Lemma \ref{L:timedependent}, all we need to show is that for any function $\phi\in\Ce_0^\infty(\R^{2d})$ on phase space,
the expected value 
\begin{equation*}
   f_{\eps,\phi}(t)=\int_{\R^{2d}}W_\eps(x,p,t)\phi(x,p)\d(x,p)
\end{equation*}
is differentiable in $t$ and satisfies hypothesis \eqref{H2} 
of Lemma \ref{L:timedependent} (ii), i.e.
$f_{\eps,\phi}^\prime(t)$ stays bounded
independently of $\eps$ and $t$. 
Differentiability in $t$ holds even for $\phi\in\L^1(\R^{2d})$,
since $W_\eps\in\Ce^1(\R;\L^\infty(\R^{2d}))$ by Lemma \ref{L:rigorousWigner}.

To deduce \eqref{H2} we start by exploiting the Wigner equation \eqref{WE}. This yields
\begin{align*}
   |f_{\eps,\phi}^\prime(t)| \le &\,
      \Big|\int_{\R^{2d}}\Weps(x,p,t)\, p\cdot\nabla_x\phi(x,p)\d (x, p)\Big|\\
    & + \Big|\frac1{(2\pi)^{d}}
      \int_{\R^{2d}}\frac{U_b(x+\frac{\eps y}2)-U_b(x-\frac{\eps y}2)}\eps
      \Psieps(x+\frac{\eps y}2,t)\overline{\Psieps(x-\frac{\eps y}2,t)}
      (\calF_p\phi)(x,y)\d(x,y)\Big|\\
    & + \Big|\frac1{(2\pi)^{d}} \int_{\R^{2d}}\frac{U_s(x+\frac{\eps y}2)-U_s(x-\frac{\eps y}2)}\eps
      \Psieps(x+\frac{\eps y}2,t)\overline{\Psieps(x-\frac{\eps y}2,t)}
      (\calF_p\phi)(x,y)\d(x,y)\Big|.
\end{align*}
The first term on the right hand side is bounded by
$\frac1{(2\pi)^d}\|\Psieps(\cdot,t)\|^2 \| p{\cdot}\nabla_x\phi \|_{{\cal A}}$
by (\ref{elest}), and hence bounded independently of $\epsilon$ and $t$, 
{  by \eqref{estSchroed} and $\|\Psieps^0\|=1$.}

Thanks to the elementary inequality 
$|\frac{U_b(x+\frac{\eps y}2)-U_b(x-\frac{\eps y}2}{\eps}|
\le \|\nabla U_b\|_\infty |y|$, 
the second term is bounded independently of $\eps$ and $t$ by 
{  \begin{equation*}
   \frac1{(2\pi)^{d}}\|\nabla U_b\|_\infty
   \int_{\R^d}|y|\sup_{x\in\R^d}|(\calF_p\phi)(x,y)|\d y.
\end{equation*}

%
%
%
%
Finally, in order to estimate the third term we observe
that for $R, Q\in\R^n$
\begin{equation*}
  \frac1\eps
  \left|\frac1{\left|R+\frac{\eps Q}2\right|}
  -\frac1{\left|R-\frac{\eps Q}2\right|}
  \right|
  \le 
  \frac{
  \left|Q\right|
  }{\left|R+\frac{\eps Q}2\right|
  \left|R -\frac{\eps Q}2\right|}.
\end{equation*}
Hence, with $x=(R_1,\ldots,R_M)$, $y=(Q_1,\ldots,Q_M)$
and setting $R=R_\alpha-R_\beta,\ Q=Q_\alpha-Q_\beta\in\R^3$, 
the third term is bounded by 
\begin{align*}
  &\frac1{(2\pi)^d}\sum_{1\le\alpha<\beta\le M}C_{\alpha\beta}
   \|\frac1{|R_\alpha-R_\beta|}\Psieps(\cdot,t)\|^2
   \int_{\R^d}|Q_\alpha{-}Q_\beta|\sup_{x\in\R}|(\calF_p\phi)(x,y)|\d y\\
  &\le\frac1{(2\pi)^d} \frac{2}{m}
     \|U_s\Psieps(\cdot,t)\|^2\int_{\R^d}|y|\sup_{x\in\R}|(\calF_p\phi)(x,y)|\d y,
\end{align*}
where $m=\min\{C_{\alpha\beta} \, | \, C_{\alpha\beta}\neq 0\}$. 
But the right hand side stays bounded independently of $\eps$ and $t$ thanks to 
Lemma \ref{bddUpsi}, establishing hypothesis \eqref{H2} 
and thus completing the proof that 
$W\in \Ce_{weak*}(\R;\calM(\R^{2d}))$, $W(t)\ge 0$.
\\[2mm]
It remains to show that $\int_{\R^{2d}}\d W(t)=1$ for every $t$. 
{In Lemma \ref{L:timedependent} (ii) we proved}
 that $W_\eps(\cdot,t)$ converges weak* in $\calA'$ to $W(t)$ for every $t$. 
Hence
Lemma \ref{L:stationary} (iv) is applicable and it suffices to verify that
both the sequence of position densities $|\Psi_\eps(\cdot, t)|^2$ 
and momentum densities 
$|\frac1{(2\pi\eps)^{d/2}}(\calF\Psi_\eps)(\frac{p}\eps,t)|^2$ are tight. 
As regards the momentum densities, this follows from 
uniform boundedness of kinetic energy, (\ref{est5}), since
\begin{align*}
 \frac{1}{2}\int_{|p|\ge R} 
\Bigl|\frac1{(2\pi\eps)^{d/2}}(\calF\Psi_\eps)(\textstyle{\frac{p}{\eps}},t)
\Bigr|^2\d p 
 &  \le \frac{1}{R^2} \frac{1}{2}\int_{\R^d} 
\Bigl|\frac1{(2\pi\eps)^{d/2}}(\calF\Psi_\eps)(\textstyle{\frac{p}{\eps}},t)
\Bigr|^2 {|p|^2} \d p  \\
 &   =  \frac{1}{R^2} \frac{1}{2}\int_{\R^d} |\eps\nabla\Psi_\eps(x,t)|^2 \d x.
\end{align*}
Finally, tightness of the position densities follows (under much weaker hypotheses on potential
and initial data) from the Lemma \ref{L:tightness} below, 
completing the proof of Theorem \ref{limitequation} (ii).\hfill$\square$
\begin{lemma}\label{L:tightness} (Propagation of tightness) 
Let $U$ be {  as in \eqref{katocdn},} 
and suppose in addition
$U$ bounded from below. Let $\{\Psieps\}$ be a sequence of solutions
to the time-dependent Schr\"odinger equation (\ref{SE}), whose initial data satisfy
$\Psi_\eps^0\in \H^2(\R^d)$, $\|\Psi_\eps^0\|=1$ (normalization), and 
$\langle\Psi_\eps^0, \, \Heps\Psi_\eps^0\rangle \le C$ (bounded energy).  
If the sequence of initial position densities $\{|\Psieps^0|^2\}$ is tight, 
then so is $\{|\Psieps(\cdot, t)|^2\}$, for all $t\in\R$. 
\end{lemma}
{\bf Proof of Lemma \ref{L:tightness}}
Let $\chi\in \Ce_0^\infty(\R^d)$, $0\le\chi\le 1$, $\chi=1$ on $|x|>1$, $\chi=0$ on $|x|<1/2$, and for $R>0$
set $\chi_R(x):=\chi(x/R)$. Then $|\nabla\chi_R|\le C/R$, 
$|\Delta\chi_R|\le C/R^2$ for some constant $C$ independent of $R$. 
Letting $\langle\cdot , \cdot \rangle$ be the 
$\L^2(\R^d)$ inner product and abbreviating
$\langle A\rangle_{\phi}=\langle\phi, \, A\phi\rangle$,
we obtain from \eqref{SE}
\begin{equation*}
\frac\dd{\dd t}\langle\chi_R\rangle_{\Psieps(\cdot, t)}
=\frac{i}{\eps}\langle[\chi_R,\Heps]\rangle_{\Psieps(\cdot, t)}.
\end{equation*}
Since 
$[\chi_R,\Heps]=\frac{\eps^2}{2}\Delta\chi_R+\eps^2\nabla\chi_R\cdot\nabla$, 
\begin{align*}
\frac\dd{\dd t}\langle\chi_R\rangle_{\Psieps(\cdot,t)}
&
\le \frac\eps2\int_{\R^d}|\Delta \chi_R(x)|\,|\Psieps(x,t)|^2\d x
+\int_{\R^d}|\eps\nabla\Psieps(x,t)|\,|\nabla \chi_R(x)|\,|\Psieps(x,t)|\d x
\\&
\le \frac\eps2\,\|\Delta \chi_R\|_\infty\,\|\Psieps(\cdot,t)\|
+\|\nabla \chi_R\|_\infty\,\|\eps\nabla\Psieps(\cdot,t)\|\,
\|\Psieps(\cdot,t)\|. 
\end{align*}
>From the boundedness from below
of $U$ and the conservation in time of the energy 
$\langle\Psi_\eps(\cdot,t), H_\eps\Psi_\eps(\cdot, t)\rangle$,
we obtain $\|\eps\nabla\Psieps(\cdot,t)\|\le const.$
for some constant independent of $\eps$ and $t$. 
Using the bounds on $\|\nabla\chi_R\|_\infty$ and $\|\Delta\chi_R\|_\infty$,
and assuming without loss of generality $\eps\le R$,}
it follows that $\frac{\dd}{\dd t}\langle \chi_R\rangle_{\Psieps(\cdot, t)}
\le\frac{const.}R$ for some constant independent of $\eps$ and $t$. 
Consequently (considering without
loss of generality $t\ge 0$) 
\begin{align*}
    \int_{|x|>R}|\Psieps(x,t)|^2\d x 
    & \le 
\langle \chi_R\rangle_{\Psieps(\cdot, t)}
\le\langle \chi_R\rangle_{\Psieps^0}
+ \int_{s=0}^t \frac{const.}{R} \d s 
\le\int_{|x|>\frac{R}2}|\Psieps^0(x)|^2\d x 
+ \frac{const.}{R}t 
\to 0 
\end{align*}
as $R\to\infty$,
by the tightness of the sequence  of initial position densitites.\hfill$\square$ 
%
%
%
%
\section{Justification of the Liouville equation}
We now show that time-dependent Wigner measures arising as a limit of semiclassically scaled
solutions to the Schr\"odinger equation are weak solutions of the Liouville equation \eqref{LE}. 

In line with the discussion in the Introduction, the main work will go into 
analyzing the behaviour near the Coulomb singularities. The first order of business will be to
verify that the limit equation even makes sense, for which we need that the gradient of the potential
lies in $\L^1$ with respect to the Wigner measure, for all $t$. 

We note the particularly impeding feature
that the singularities of the Coulomb forces appearing in the Liouville equation
have the magnitude $1/|R_\alpha-R_\beta|^2$, and are hence much worse than the
singularities of the Coulomb potentials $1/|R_\alpha-R_\beta|$ appearing in the original Schr\"odinger equation.
\\[2mm]
{\bf Proof of Theorem \ref{limitequation} (iii)} We want to show that 
\begin{eqnarray}
 & & W(t)({\cal S}\times\R^d)=0, \label{null} \\
 & & \int_{\R^{2d}} |\nabla U_s| \d W(t) < \infty \quad\text{for all $t\in\R$.} \label{L1}
\end{eqnarray}
To this end, we fix $t\in\R$ and choose a subsequence such that $|\Psieps(\cdot, t)|^2\weakstar\mu$ in $\calM(\R^d)$. 
{  According to Lemma \ref{L:stationary} (iii)} we have
\begin{equation}\label{mux}
   \int_{p\in\R^d}W(\cdot,\dd p, t)\le\mu.
\end{equation}
First we show $W(t)({\cal S}\times\R^d)=0$. By (\ref{mux}) it suffices to show that $\mu({\cal S})=0$. 
Let {  $\phi\in \Ce(\R^d)$} with $0\le \phi\le 1$, $\phi=1$ for 
$\dist(x,\calS)<\delta$, and $\phi=0$ for $\dist(x,\calS)\ge 2\delta$. We use the elementary property of
weak* convergence in $\calM(\R^d)$ that if $\mu_\eps\ge 0$, $\mu_\eps\weakstar\mu$ in $\calM(\R^d)$, 
and $f$ is a bounded nonnegative continuous function, then $\int_{\R^d}f \, \d\mu \le 
\liminf\limits_{\eps\to 0}\int_{\R^d}f\, \d\mu_\eps$.
(This follows immediately by approximating $f$ by compactly supported functions $f_j$ with $0\le f_j\le f$ and
dominated convergence.) Consequently
$$
   \int_{\dist(x,\calS)\le\delta}\d\mu
   \le \int_{\R^d} \phi\, \dd\mu \\
   \le {  \liminf_{\eps\to 0} \int_{\R^d} \phi\, |\Psieps(\cdot, t)|^2\d x
   =\liminf_{\eps\to 0}\int_{\dist(x,\calS)\le 2\delta}|\Psieps(\cdot,t)|^2}
   \d x.
$$
The idea now is to use the fact that $U_s\ge\frac1{\wt C\delta}$ on $\{x\in\R^d\, | \, \dist(x,\calS)\le 2\delta\}$
for some constant $\wt C$, and appeal to Lemma \ref{bddUpsi} below. This yields
\be \label{concest}
  \int_{\dist(x,\calS)\le\delta}\dd\mu \le 
  \Bigl(\wt C\delta\Bigr)^2 { 
    \liminf_{\eps\to0}}\int_{\R^d}U_s^2|\Psieps(\cdot, t)|^2 \d x\le 
  \Bigl(\wt C\delta\Bigr)^2 C.
\ee
Since the RHS tends to zero as $\delta\to 0$, $\mu(\calS)=0$ and $W(t)(\calS\times\R^d)=0$, completing the
proof of (\ref{null}). 

To show (\ref{L1}) we will use the monotone convergence theorem. 
To this end we set for $\delta>0$ 
\begin{equation*}
   f_\delta:=\min\Bigl\{f,{\frac1\delta}\Bigr\},\quad f:=|\nabla U_s|.
\end{equation*}
Then $f_\delta$ is a bounded continuous function on $\R^d$. 
By the fact that $|\nabla U_s|\le C_0 U_s^2$ (thanks to the identity
$\left|\nabla_{R_\alpha}\frac1{|R_\alpha-R_\beta|}\right|
=\frac1{|R_\alpha-R_\beta|^2}$) and Lemma \ref{bddUpsi} below, 
\begin{equation} \label{mux2}
   \int_{\R^d}f_\delta \d\mu \le
   \liminf_{\eps\to 0}\int_{\R^d}f_\delta|\Psieps(\cdot, t)|^2\d x 
   \le C_0\liminf_{\eps\to 0}\int_{\R^d}U_s^2|\Psieps(\cdot, t)|^2
   \d x\le C_0 C,
\end{equation}
where $C_0$, $C$ are constants independent of $\eps$ and $t$.

Consider now the limit $\delta\to 0$. In this limit, $f_\delta\to f$ monotonically on $\R^d\setminus\calS$.
Since $\mu(\calS)=0$, it follows that $f_\delta\to f$ $\mu$-almost everywhere.
Hence the monotone convergence theorem yields $f\in\L^1(\d\mu)$ and
\be \label{mux3}
   \int_{\R^d} f\d\mu=\lim_{\delta\to0}\int_{\R^d} f_\delta\d\mu.
\ee
Consequently, by the $p$-independence of $f$, 
(\ref{mux}), (\ref{mux3}), (\ref{mux2}),  
\be \label{fbound}
   \int_{\R^{2d}} f  \d W(t)
   \le \int_{\R^{d}}f\d\mu \le C_0 C <\infty,
\ee
establishing (\ref{L1}) and completing the proof of (iii).\hfill$\square$
\\[2mm]
{\bf Proof of Theorem \ref{limitequation} (iv)} 
Let $\Omega\subseteq\R^d$ be an open set such that $U_b\in\Ce^1(\Omega)$. We need to show that eq. (\ref{weakLE})
holds for every $\phi\in\Ce_0^\infty(\Omega{\times}\R^d{\times}\R)$. Starting point is the fact that
by Lemma \ref{L:rigorousWigner}, $\Weps$ satisfies the Wigner equation (\ref{WE}). 
Multiplying (\ref{WE}) by $\phi$ and integrating by parts, we obtain the weak form 
\be \label{weakWE}
    \int_{\R^{2d+1}}
    \Bigl(\Weps(\partial_t + p \cdot \nabla_x)\phi + f_\eps \phi\Bigr)\d (x,p,t) = 0.
\ee
Passage to the limit $\eps\to 0$ in eq. \eqref{weakWE} is done in three steps, carried out
in the order of increasing difficulty: 1. Analysis of the
local terms, 2. Analysis of the nonlocal term $\int f_\eps\phi$ for test functions vanishing in a
neighbourhood of the Coulomb singularities, 3. Analysis of the nonlocal term in a neighbourhood of the
Coulomb singularities. 
\\[2mm]
{\bf Step 1: Analysis of the local terms in \eqref{weakWE}} \\[1mm]
As $\eps\to 0$, $W_\eps\weakstar W$ in $\L^\infty(\R;\calA^\prime)$,
and hence in particular in ${\cal D}'(\R^{2d+1})$.
Thus, the first term in \eqref{weakWE} satisfies 
\begin{equation}\label{wignerlimit1}
    \int_{\R^{2d+1}}
    \Weps\,(\partial_t + p \cdot \nabla_x)\phi \d (x,p,t)
\to \int_{\R}\int_{\R^{2d}}
    (\partial_t + p \cdot \nabla_x)\phi \d W(t)\d t
\quad(\eps\to 0).
\end{equation}
{\bf Step 2: Analysis of the nonlocal term $\int f_\eps\phi$ for test functions vanishing in a
neighbourhood of the Coulomb singularities} \\[1mm]
Let ${\cal S}\subset\R^d$ be the set of Coulomb singularities
of the potential $U_s$ (see Theorem \ref{limitequation} (iii)). Our goal in this step is to prove
that 
\be \label{wignerlimit2} 
   \int_{\R^{2d+1}} f_\eps\,\phi\d (x,p,t) \to - \int_{\R}\int_{\R^{2d}} 
\nabla U \cdot \nabla_p\phi \d W(t)\d t 
\quad(\eps\to 0)
\ee
for all test functions on $\Omega\times\R^d\times\R$ 
which vanish in a neighbourhood of $\calS$, i.e., 
$\phi\in\Ce_0^\infty((\Omega\backslash{\cal S})\times\R^d\times\R)$. 
This together with (\ref{weakWE}), 
(\ref{wignerlimit1}) completes the proof of (iv) for the above test functions. 

We begin by rewriting the left hand side of (\ref{wignerlimit2}). Substituting the definition 
of $f_\eps$ and carrying out the integration over $p$ gives
(abbreviating $\Psieps = \Psieps(x+\frac{\eps y}2,t)$, $\overline{\Psieps}=
\overline{\Psieps(x-\frac{\eps y}2,t)}$, $\calF_p\phi = (\calF_p\phi)(x,y,t)$)
\be \label{reformul}
   \int_{\R^{2d+1}}f_\eps\,\phi \d (x,p,t)= 
   -\frac{i}{(2\pi)^d}\int_{\R^{2d+1}} \!\!
   \frac{U(x+\frac{\eps y}2) - U(x-\frac{\eps y}2)}{\eps} 
   \Psieps\overline{\Psieps} 
   (\calF_p\phi)\d (x,y,t).
\ee

The idea now is to split the $y$-integration into two regions in such a way that either 
$|\eps y|<<1$, in which case the difference quotient of $U$ is well approximated by the derivative
$\nabla U(x)\cdot y$, or $|y|>>1$, in which case $\calF_p\phi$ is very small, due to the rapid decay
of the Schwartz function $\calF_p\phi$ as $|y|\to\infty$.
To implement this idea, we introduce a cut off radius 
$\eps^{-\alpha}$ with some fixed $\alpha\in(0,1)$  
and choose the regions of integration as $|y|\le \eps^{-\alpha}$ and 
$|y|>\eps^{-\alpha}$. 
{  
In particular, we denote 
$\Lambda^\eps:=\{(y,t)\in\R^{d+1}\,|\,|y|\le\eps^{-\alpha}\}$
and $V_\eps:=\{(y,t)\in\R^{d+1}\,|\, |y|>\eps^{-\alpha}\}$.
}
It is also convenient to subtract off, and then add again,
the RHS of \eqref{reformul} with the difference quotient of $U$ replaced by the derivative
$\nabla U(x)\cdot y$. This yields the following natural splitting of the RHS of \eqref{reformul}
into a sum of five terms 
{  
\begin{align} 
\label{splitting}
  &\int_{\R^{2d+1}}  f_\eps\, \phi  \d (x,p,t) 
=\  T_1 + T_2^+ + T_2^- + T_3 + T_4, \\
& T_1 := 
-\frac{i}{(2\pi)^d}
\int_{\R^d\times \Lambda^\eps}
\Bigl[ \frac{U(x+\frac{\eps y}2) - U(x-\frac{\eps y}2)}{\eps} - \nabla U(x) 
\cdot y\Bigr]\frac{1}{|y|} \cdot |y|\,  
   \Psieps\overline{\Psieps} (\calF_p\phi) \d (x,y,t), \nonumber \\
& T_2^{\pm} := \mp \frac{i}{(2\pi)^d}
\int_{\R^d\times V_\eps}
\frac{U(x\pm\frac{\eps y}2)}{\eps}\,  \Psieps\overline{\Psieps} 
 (\calF_p\phi)  \d (x,y,t), \nonumber \\
 & T_3 :=  \frac{i}{(2\pi)^d}
\int_{\R^d\times V_\eps}
     \nabla U(x)\cdot y \, \Psieps\overline{\Psieps} 
   (\calF_p\phi)  \d (x,y,t), \nonumber \\
  & T_4 := - \frac{i}{(2\pi)^d}\int_{\R^{2d+1}} \nabla U(x)\cdot y \, \Psieps\overline{\Psieps} 
   (\calF_p\phi)  \d (x,y,t). \nonumber 
\end{align}
}
Here, in $T_1$ the factor $\frac{1}{|y|}\cdot |y|=1$ has been inserted for future use. 

We begin by analysing $T_1$. 
Let $\delta, \, R>0$ such that the set 
$\Omega_{R,\delta} :=\{ x\in\Omega \, | \, |x|\le R, \, 
\dist(x,\calS\cup\partial\Omega)\ge\delta\}$ 
contains all $x$ with $(x,p,t)\in\supp\phi$.
Then provided $\eps$ is sufficiently small
$$
   |x\pm\frac{\eps y}2|\le 2R, \;\;  \dist(x\pm\frac{\eps y}2,\calS\cup\partial\Omega)\ge\frac\delta2 
   \quad\mbox{ for all }x\in \Omega_{R,\delta}, \, |y|\le \eps^{-\alpha},
$$
that is to say $x\pm\frac{\eps y}2\in\Omega_{2R,\delta/2}$. 
Hence by the continuous differentiability of $U$ in $\Omega_{2R,\delta/2}$, 
$$
   \frac{U(x+\frac{\eps y}2) - U(x-\frac{\eps y}2)}{\eps|y|}
   \to \nabla U(x)\cdot \frac{y}{|y|}
   \quad\text{ as $\eps\to0$ uniformly for }x\in\Omega_{R,\delta}, \, |y|\le\eps^{-\alpha}.
$$
Consequently, applying 
the Cauchy-Schwarz inequality with respect to the integration over $x$, 
$T_1$  can be estimated by
\be \label{wignerlimit4}
  |T_1| \le \frac1{(2\pi)^d} 
   \sup_{x\in\Omega_{R,\delta}, \, |y|\le \eps^{-\alpha}} \!
   \left|\frac{U(x+\frac{\eps y}2) - U(x-\frac{\eps y}2)}{\eps|y|} 
   -\nabla U(x)\cdot\frac{y}{|y|}\right| 
\,\sup_{t\in\R}\|\Psieps(\cdot,t)\|^2 
\,\|y \calF_p\phi \|_\ast
   \to 0 
\ee
as $\eps\to0$, 
where here and below, for any function $\chi\in\calS(\R^{2d+1})$ we denote
\begin{equation}\label{normast}
   \|\chi\|_\ast:=\int_{\R^{d+1}}\sup_{x\in\R^d}|\chi(x,y,t)|\d(y,t).
\end{equation}
The terms $T_2^\pm$ can be estimated in an analogous manner, 
again applying the Cauchy-Schwarz inequality
with respect to the integration over $x$: 
\be \label{Ttwoest1}
  |T_2^\pm| \le 
\frac1{(2\pi)^{d}}
\,\sup_{t\in\R}\|U\Psieps(\cdot, t)\|\,\sup_{t\in\R}\|\Psieps(\cdot, t)\|
\,\frac1{\eps}
\int_{V_\eps}
\sup_{x\in\R^d}|(\calF_p\phi)(x,y,t)|
\d(y,t). 
\ee
By Lemma \ref{bddUpsi} and the boundedness of $U_b$, the norm 
$\|U\Psieps(\cdot, t)\|$ stays bounded independently of $t$ and $\eps$.
%
On the other hand,
since $\calF\phi \in {\cal S}(\R^{2d+1})$, 
$$\displaystyle 
   \sup_{(x,y,t)\in\R^{2d+1}} (1{+}|y|)^m |(\calF_p\phi)(x,y,t)| =: c_m <\infty
$$
for any $m=0,1,2,\ldots$. Consequently, for all $m\ge d+1$ 
and all $|y|>\eps^{-\alpha}$
$$
   |(\calF_p\phi)(x,y,t)| \le \frac{c_m}{(1{+}|y|)^m} 
   \le \frac{c_m}{(1{+}|y|)^{d+1}} \frac{1}{(\eps^{-\alpha})^{m-d-1}} 
= \frac{c_m}{(1{+}|y|)^{d+1}} \eps^{\alpha(m-d-1)}.
$$
Hence, choosing $T$ so large that $|t|\le T$ for all 
$(x,p,t)\in\supp\phi$, the last factor in \eqref{Ttwoest1} satisfies 
{  
\be \label{tailest}
\int_{V_\eps}
\sup_{x\in\R^d}|(\calF_p\phi)(x,y,t)|
\d(y,t) 
  \le 
\eps^{\alpha(m-d-1)}\int_{\R^d}\frac{2 T c_m}{(1{+}|y|)^{d+1}}\d y.
\ee
}
Then, if $m$ is chosen 
so large that the exponent $\alpha(m-d-1)>1$,
the positive power of $\eps$ in \eqref{tailest} 'beats' the singular factor
$\frac1\eps$ in \eqref{Ttwoest1}.
Thus 
\be \label{wignerlimit5}
   T_2^\pm \to 0 \mbox{ as }\eps\to 0
\ee
(in fact, faster than any power of $\epsilon$, but this is not needed in the sequel). 

The proof that $T_3\to 0$ is analogous but easier, due to the absence of a 
singular prefactor $\frac{1}{\eps}$ and the fact that $\nabla U(x)$, 
unlike $U(x\pm\eps y/2)$, is bounded on $\supp\phi$. We simply estimate
\be \label{wignerlimit6}
  |T_3| \le \frac{1}{(2\pi)^d} \sup_{x\in\Omega_{R,\delta}}|\nabla U(x)| 
\int_{V_\eps}
|y|\sup_{x\in\R^d}|(\calF_p\phi)(x,y,t)|
\d(y,t)\to 0\mbox{ as }\eps\to 0
\ee
by \eqref{tailest} {  with $m>d+2$.} 

Finally, consider the last term, $T_4$. Using 
$iy \calF_p\phi=\calF_p(\nabla_p\phi)$ and 
interchanging the
integrations over $p$ and $y$,
\begin{align*}
 T_4 =& -\frac{1}{(2\pi)^d}\int_{\R^{2d+1}}\nabla U(x)\cdot \Bigl[
 \int_{\R^d} \nabla_p \phi(x,p,t) e^{-ip\cdot y} \d p
 \Bigr] \Psieps(x+\frac{\eps y}2,t)\overline{\Psieps(x-\frac{\eps y}2,t)}
 \d (x,y,t) \\
     =& -\int_{\R^{2d+1}} \nabla U(x)\cdot \nabla_p \phi(x,p,t)
         \, \Weps(x,p,t) \d (x,p,t).
\end{align*}
Since $\nabla U \cdot \nabla_p \phi \in \L^1(\R;\calA)$ due to the continuity 
of  $\nabla U$ on $\supp\phi$, and since $\Weps\weakstar W$ in the dual 
$\L^\infty(\R;\calA^\prime)$,
\be \label{wignerlimit7}
  T_4 \to -\int_{\R}\int_{\R^{2d}}\nabla U(x)\cdot \nabla_p \phi(x,p,t)
 \d W (t) \d t \mbox{ as $\eps\to0$.}
\ee
Combining \eqref{wignerlimit4}, \eqref{wignerlimit5}, \eqref{wignerlimit6}, \eqref{wignerlimit7} yields (\ref{wignerlimit2}) for 
$\phi\in\Ce_0^\infty((\Omega\setminus\calS){\times}\R^d{\times}\R)$.
\\[2mm]
{\bf Step 3: Analysis of the nonlocal term $\int f_\eps\phi$ in a
neighbourhood of the Coulomb singularities} \\[1mm]
We prove here that eq. (\ref{wignerlimit2}) continues to hold for arbitrary 
$\phi\in\Ce_0^\infty(\Omega\times\R^d\times\R)$ 
not required to vanish in a neighbourhood of the Coulomb singularities, 
i.e. that the Liouville equation continues to hold across
Coulomb singularities. 

This is quite remarkable, since the available a priori bound
$\int_{\R^d}U_s^2|\Psieps(\cdot, t)|^2\le const.$ 
only rules out concentration of the measure $|\Psieps(\cdot, t)|^2$ 
on Coulomb singularities (as was shown in (\ref{concest})), 
but not concentration of the blown-up measure 
$U_s^2|\Psieps(\cdot, t)|^2\sim |\nabla U_s| \, |\Psieps(\cdot, t)|^2$ 
which asymptotically appears in $\int f_\eps\phi$ 
(see the leading term $T_4$ in (\ref{splitting})). 
This suggests the possibility that an additional contribution of form 
$\int_\R\int_{\calS\times\R^d}\phi  \d\nu(t) \d t$, 
with
$\nu(t)$ a singular measure supported on $\calS\times\R^d$, could appear in the limit equation 
(\ref{weakLE}). 
 
This possibility will be ruled out by careful use of the evolution equation (\ref{WE}) satisfied
by $W_\eps$. Roughly, our analysis below will lead  
to the insight that the asymptotic amount of mass of $f_\eps$ in a $\delta$-neighbourhood of
the singular set $\calS$ is at most of order $\delta$, not order one. 

We will need rather precisely chosen cutoff functions. Given $\delta>0$, we let
{ 
$$
\eta_\delta(x) 
:= \prod_{1\le\alpha<\beta\le M} \eta\Bigl(\frac{|R_\alpha-R_\beta|}{\delta}
\Bigr)
\quad \text{with $x=(R_1,\ldots,R_m)$,}
$$
where $\eta\in \Ce_0^\infty(\R)$ with $0\le\eta\le 1$, $\eta=1$ on 
$|z|\le 1/2$, $\eta=0$ on $|z|\ge 1$. 
} 
Then for some constants $C_1$, $C_2$ independent of $\delta$
\be \label{etabound}
   |\nabla\eta_\delta(x)|\le \frac{C_1}{\delta} \ \mbox{ for all } x, 
\quad \eta_\delta(x)=0\ \mbox{ for }\dist(x,\calS)\ge C_2\delta.
\ee
We now write $\phi=(1-\eta_\delta)\phi + \eta_\delta\phi$, and consider 
both contributions to $\int_{\R^{2d+1}}f_\eps\phi \d(x,p,t)$ separately. 
Since $(1-\eta_\delta)\phi$ belongs to 
$\Ce_0^\infty((\Omega\setminus\calS)\times\R^d\times\R)$, 
we have by Step 2 (cf.\ \eqref{wignerlimit2})
\be \label{wignerlimit8}
   \int_{\R^{2d+1}} f_\eps \, (1-\eta_\delta)\phi \d (x,p,t) \to 
   - \int_{\R}\int_{\R^{2d}}\nabla U \cdot (1-\eta_\delta)\nabla_p\phi \d
   W(t)\d t
\quad\text{as $\eps\to 0$.} 
\ee
We now claim that the remaining terms are small when $\delta$ is small, that is to say
\begin{eqnarray} \label{*}
  & & \limsup_{\eps\to 0} \Bigl|\int_{\R^{2d+1}}f_\eps\, \eta_\delta\phi 
\d (x,p,t) \Bigr| \to 0 \mbox{ as }\delta \to 0, \\
  & & \int_\R \int_{\R^{2d}} \nabla U\cdot  \eta_\delta\nabla_p\phi 
\d W(t)\d t \to 0
  \mbox{ as }\delta\to 0. \label{*2}
\end{eqnarray}
Clearly, (\ref{wignerlimit8}) and (\ref{*}), (\ref{*2}) imply eq. (\ref{wignerlimit2}) for the arbitrary test function
$\phi\in\Ce_0^\infty(\Omega\times\R^d\times\R)$. 

Proving (\ref{*2}) is not difficult, but for convenience of the reader
we include a proof. Consider first a fixed $t$. 
By  the facts that $0\le\eta_\delta\le 1$ and
$\supp\eta_\delta\subset\{x\in\R^d\,|\,\dist(x,\calS)\le C_2\delta \}$,
\begin{align} \label{fest2}
\Bigl|\int_{\R^{2d}}\nabla U\cdot\eta_\delta\nabla_p\phi\d W(t)\Bigl|
\,\le\, & \|\nabla_p\phi\|_\infty 
\int_{(\supp\eta_\delta\cap\supp\phi)\times\R^d}
|\nabla U |\d W(t)
 \\& 
\to\|\nabla_p\phi\|_\infty
\int_{(\calS\cap\supp\phi)\times\R^d}|\nabla U |\d W(t) = 0 
\mbox{ as }\delta\to 0, \nonumber
\end{align}
since $\d W(t)(\calS\times\R^d)=0$ and $\nabla U\in \L^1(\supp\phi;\d W(t))$. 
%
Moreover the LHS of (\ref{fest2}) stays bounded independently of $t$ 
by (\ref{fbound}). 
Hence the integrand with respect to $t$
in (\ref{*2}) tends to zero boundedly a.e.\ as $\delta\to 0$, and so by dominated convergence 
we infer (\ref{*2}).

It remains to establish (\ref{*}). 
This is the difficult part of Step 3, due to the fact discussed above
that the a priori bound of Lemma \ref{bddUpsi} does not rule out 
the possibility of concentration of mass of $|\nabla U\|\Psi_\eps(\cdot,t)|^2$ 
on the set $\calS$ of Coulomb singularities. 
Using first (\ref{WE}) and then the definition of $W_\eps$ we have
{ 
\begin{align*}
& \int_{\R^{2d+1}}f_\eps\, \eta_\delta\phi \d(x,p,t) 
= -\int_{\R^{2d+1}} \Weps \,(\partial_t + p\cdot
  \nabla_x)(\eta_\delta\phi)\d(x,p,t) 
\\ & 
=  -\int_{\R^{2d+1}}  W_\eps 
\Bigl[\eta_\delta(\partial_t+p\cdot\nabla_x)\phi + \nabla\eta_\delta \cdot
                           p\phi\Bigr] \d(x,p,t) 
=  -\frac{1}{(2\pi)^d} \int_{\R^{2d+1}} u_{\eps,\delta}\d(x,y,t)
\\& 
= \underbrace{-\frac{1}{(2\pi)^d} 
\int_{\R^d\times\Lambda^\eps} 
u_{\eps,\delta}\d(x,y,t)}_{=:Q_1} 
\underbrace{-\frac{1}{(2\pi)^d} 
\int_{\R^d\times V_\eps}
u_{\eps,\delta} \d(x,y,t)}_{=:Q_2}, 
\end{align*}
}
where we have split the domain of integration as in Step 2, and where
$$
   u_{\eps,\delta}(x,y,t) = \Bigl[ \eta_\delta(x) 
\underbrace{\calF_p((\partial_t+p\cdot\nabla_x)\phi)}_{=:\chi}(x,y,t) 
                +\nabla\eta_\delta(x)\cdot 
\underbrace{\calF_p(p\phi)}_{=:\xi}(x,y,t)\Bigr] 
\Psieps(x+\frac{\eps y}{2},t)
                 \overline{\Psieps(x-\frac{\eps y}{2},t)}.
$$

The second term, $Q_2$, is exactly of the same form as the term $T_3$ in Step 2, with $\nabla U(x)$ replaced by 
$\eta_\delta(x)$ respectively $\nabla\eta_\delta(x)$, 
and the test function $y\calF_p\phi$ replaced by 
$\chi$ respectively $\xi$. Consequently, estimating
as for (\ref{wignerlimit6}), 
{  
\begin{equation}
|Q_2| \le \frac{1}{(2\pi)^d} 
\Bigl[\|\eta_\delta\|_\infty 
\int_{V_\epsilon}\sup_{x\in\R^d}|\chi(x,y,t)|\d(y,t) 
+ \|\nabla\eta_\delta\|_\infty
\int_{V_\epsilon}\sup_{x\in\R^d}|\xi(x,y,t)|\d(y,t)\Bigr]  
\to  0 \; (\eps\to 0),  
\label{wignerlimit9}
\end{equation}
by (\ref{tailest}) with $\chi$ respectively $\xi$ in place of $\calF_p\phi$,
and $m>d+1$. 
}
The first term, $Q_1$, will be dealt with by an argument similar to the no-concentration estimate (\ref{concest})
in the proof of (iii), except now a non-local version is needed since the term under investigation is not local in $W_\eps$. 
For $\eps$ sufficiently small, and $(x,y)$ belonging to
the domain of integration, that is to say $x\in\supp\eta_\delta$, 
$|y|\le \eps^{-\alpha}$, we have (writing $x=(R_1,..,R_M)$, 
$y=(Q_1,..,Q_M)$) 
$$
   |(R_\alpha \pm \frac{\eps Q_\alpha}{2}) - (R_\beta \pm \frac{\eps
     Q_\beta}{2})| \le 2\delta\quad\text{for $1\le\alpha<\beta\le M$}
$$
and consequently $U_s(x \pm \frac{\eps y}{2})\ge\frac{1}{{\wt C}\delta}$ for 
$x\in\supp\eta_\delta$, $|y|\le\eps^{-\alpha}$,
and some constant $\widetilde C$.
Hence 
{  
\begin{align*}
|Q_1| &\le \frac{({\wt C}\delta)^2}{(2\pi)^d} 
\int_{\R^d\times\Lambda^\eps}
\Bigl[\eta_\delta|\chi|+|\nabla\eta_\delta|\,|\xi|\Bigr]
\,|U_s(x{+}\frac{\eps y}{2})\Psi_\eps(x{+}\frac{\eps y}{2},t)|
\,|U_s(x{-}\frac{\eps y}{2})\Psi_\eps(x{-}\frac{\eps y}{2},t)| 
\d(x,y,t)  
\\
& \le \frac{({\wt C}\delta)^2}{(2\pi)^d} 
\Bigl[\|\eta_\delta\|_\infty
\|\chi\|_\ast
+ \|\nabla\eta_\delta\|_\infty
\|\xi\|_\ast
\Bigr]\, \sup_{t\in\R} \|U_s\Psi_\eps(\cdot, t)\|^2
\end{align*}
with $\|\cdot\|_\ast$ given by \eqref{normast}.
}
Since $\|U_s\Psi_\eps(\cdot, t)\|$ stays bounded independently of $t$ and 
$\epsilon$ by Lemma \ref{bddUpsi}, 
$\|\eta_\delta\|_\infty=1$, and 
$\|\nabla\eta_\delta\|_\infty\le\frac{C_1}{\delta}$ by (\ref{etabound}), 
it follows that
\be \label{wignerlimit10}
   |Q_1|\le C_*[\delta^2 + \delta] 
\ee
for some constant $C_*$ independent of $\epsilon$ and $\delta$. Note how the positive power of $\delta$
gained by inserting the multiplier $U_s$ has `beaten' the 
negative power of $\delta$ coming
from the gradient of the cutoff function $\eta_\delta$. 
Combining (\ref{wignerlimit9}), (\ref{wignerlimit10})
gives \eqref{*}. This completes the proof of Theorem \ref{limitequation} (iv) for general test functions.\hfill$\square$ 
%
%
%
%
\section{An a priori estimate for the Schr\"odinger equation with repulsive 
Coulomb interactions} 

We prove now the a priori estimate used in the derivation of the Liouville equation that
the potential term {  $U_s\Psi_\eps(\cdot, t)$} 
in the semiclassically scaled Schr\"odinger equation (\ref{SE}) 
stays bounded in $\L^2(\R^d)$ independently of $\epsilon$ and $t$. 

For potentials with Coulomb singularities, such as (\ref{pes}), such an estimate 
says in particular that in the limit $\eps\to 0$, the wavefunction cannot concentrate mass at the singularities. 

On physical grounds, one would expect this to be true only for repulsive interactions (as present here), but
not for attractive interactions. The challenge then is to translate this physical intuition into a mathematical
argument fine enough to detect sign information. This is achieved by the positive commutator argument below,
which exploits not just the repulsivity (i.e., positivity) of $U_s$, but
also its special Coulombic nature.

\begin{lemma}\label{bddUpsi} Let $U=U_b+U_s$ be as in 
Theorem \ref{limitequation}. 
Let 
$\{\Psieps^0\}\subset\calD(\Heps)$ 
satisfy $\|\Psieps^0\|=1$ and $\|\Heps\Psieps^0\|\le c$  
for some constant $c$ independent of $\eps$.
Then the solution $\Psieps(\cdot,t)$ to \eqref{SE} satisfies
\be \label{estSchroed_add}
 \sup\limits_{t\in\R}\|U_s\Psieps(\cdot,t)\|^2\le C
\ee
for some constant $C$ independent of $\eps$.
\end{lemma}
{\bf Proof} 
By standard results on the unitary propagator $e^{-itH}$ associated to a
self-adjoint operator $H$, if $\psi_0\in\calD(H)$, then so is  
$\psi(\cdot,t)=e^{-itH}\psi_0$ for all $t\in\R$, and 
$\|\psi(\cdot,t)\|$, 
$\langle\psi(\cdot,t),H\psi(\cdot,t)\rangle$,  
$\|H\psi(\cdot,t)\|$ are time-independent. 
(Formally, the time-independence follows from the Heisenberg evolution
equation for expected values,  
$\frac{\dd}{\dd t}\langle\psi(\cdot,t),A\psi(\cdot,t)\rangle
=\langle\psi(\cdot,t),\frac1i[A,H]\psi(\cdot,t)\rangle$, 
by taking $A=I,\, H,\, H^2$.)
Applied to our case this yields, besides \eqref{estSchroed},
\begin{align}
&\sup_{t\in\R}\Bigl(\frac12\int_{\R^d}|\eps\nabla\Psieps(x,t)|^2\d x
+\int_{\R^d}(U_b(x)+U_s(x))|\Psieps(x,t)|^2\d x\Bigr)\le const.,
\label{estHpsi}\\
&\sup_{t\in\R}\Bigl\|\Bigl(-\frac{\eps^2}2\Delta+U_b+U_s\Bigr)
\Psieps(\cdot,t)\Bigr\|^2\le const.,
\label{estH2psi}
\end{align}
the constants being independent of $\eps$.
%
By \eqref{estSchroed}, the boundedness of $U_b$, and the nonnegativity 
of $U_s$, 
\begin{equation}\label{est5}
    \sup_{t\in\R}\frac12\int_{\R^d}|\eps\nabla\Psieps(x,t)|^2\d x\le const.,
 \qquad
     \sup_{t\in\R}\int_{\R^d}U_s(x)|\Psieps(x,t)|^2\d x\le const.,
\end{equation}
all constants being independent of $\eps$.
Now we expand the left hand side of \eqref{estH2psi}, and rewrite the latter in the form 
\begin{multline}\label{est4}
\sup_{t\in\R}\Big(
\Bigl\|\Bigl(-\frac{\eps^2}2\Delta+U_b\Bigr)
\Psieps(\cdot,t)\Bigr\|^2
+2\Re\Bigl\langle-\frac{\eps^2}2\Delta\Psieps(\cdot,t),U_s\Psieps(\cdot,t)\Bigr\rangle
\\
+2\Re\langle U_b\Psieps(\cdot,t),U_s\Psieps(\cdot,t)\rangle
+\|U_s\Psieps(\cdot,t)\|^2
\Big)\le const.
\end{multline}
Using the positivity of $U_s$ and \eqref{est5}, the third term satisfies
\begin{equation}\label{est6}
   \sup_{t\in\R}
   \left|2\Re\langle U_b\Psieps(\cdot,t),U_s\Psieps(\cdot,t)\rangle\right|
   \le 2\|U_b\|_\infty
\sup_{t\in\R}|\langle\Psieps(\cdot,t),U_s\Psieps(\cdot,t)\rangle|
   \le const.
\end{equation}
The key point now is the following claim:
{  
\begin{equation}\label{est7}
\Re\langle-\Delta\psi,U_s\psi\rangle\ge0\quad\text{for $\psi\in\H^2(\R^d)$.}
\end{equation}
}
Postponing its proof, 
substitution of \eqref{est6}, \eqref{est7} 
{  (with $\psi=\Psieps(\cdot,t)$)} 
into \eqref{est4} yields
$$
   \sup_{t\in\R} 
   \Bigl\|\Bigl(-\frac{\eps^2}2\Delta+U_b\Bigr)
   \Psieps(\cdot,t)\Bigr\|^2 \le const.,
\quad
   \sup_{t\in\R}\|U_s\Psieps(\cdot,t)\|^2\le const.,
$$
establishing the assertion. 

It remains to prove \eqref{est7}. (This depends on
the Coulombic nature of $U_s$ as well as the fact that it is 
positive, i.e., repulsive.)
By a standard approximation argument, using the density of $\Ce_0^\infty(\R^d)$
in $\H^2(\R^d)$ and the fact that, by Hardy's inequality, 
$\psi\mapsto U_s\psi$ is a continuous map from $\H^1(\R^d)$ to  $\L^2(\R^d)$, 
it suffices to prove \eqref{est7} for $\psi\in\Ce_0^\infty(\R^d)$.
In this case, compute
$$
  \Re\langle-\Delta\psi,U_s\psi\rangle
  =\Re\int_{\R^d}\nabla\overline{\psi}\cdot\nabla(U_s\psi)
  =\int_{\R^d}|\nabla\psi|^2 U_s
  +\Re\int_{\R^d}(\nabla\overline{\psi})\psi\cdot\nabla U_s.
$$
The first term is $\ge 0$ and the second term equals
\begin{align*}
  \frac12 \int_{\R^d}
  \bigl((\nabla\overline{\psi})\psi+\overline\psi\nabla{\psi}\bigr)
  \cdot\nabla U_s
 &=\frac12 \int_{\R^d}\nabla|\psi|^2\cdot\nabla U_s
 =\frac12 \int_{\R^d}|\psi|^2(-\Delta U_s).
\end{align*}
Now, considering e.g. the term $\ds\frac1{|R_1-R_2|}$
in $U_s$, cf.\ \eqref{Us}, 
write
\begin{equation*}
-\Delta_{(R_1,\ldots,R_M)}=-\Delta_{\frac{R_1+R_2}{\sqrt2}}
-\Delta_{\frac{R_1-R_2}{\sqrt2}}-\Delta_{(R_3,\ldots,R_M)}
\end{equation*}
and hence 
\begin{equation*}
-\Delta_{(R_1,\ldots,R_M)}\frac1{|R_1-R_2|}
=-2\Delta_{R_1-R_2}\frac1{|R_1-R_2|}
=8\pi\delta(R_1-R_2)
\end{equation*}
(recall $-\Delta\frac1{|\cdot|}=4\pi\delta$ in $\R^3$).
Thus $\int_{\R^d}|\psi|^2(-\Delta U_s)\ge 0$, 
which completes the proof of \eqref{est7}.\hfill$\square$\\
\\[4mm]
{\bf Acknowledgements} We would like to thank Caroline Lasser and Clotilde Fermanian-Kammerer for helpful discussions.

\end{document}